\documentclass{elsarticle}
\usepackage{amsmath,amssymb,amsfonts, color, inputenc, verbatim, url, graphicx, bbm}
\usepackage[skip=2pt]{caption}

\newtheorem{theorem}{Theorem}[section]
\newtheorem{prop}[theorem]{Proposition}
\newtheorem{lemma}[theorem]{Lemma}
\newtheorem{remark}[theorem]{Remark}

\newtheorem{example}[theorem]{Example}
\newtheorem{corollary}[theorem]{Corollary}

\def\cD{\mathcal{D}}

\def\cF{\mathcal{F}}

\def\cR{\mathcal{R}}

\def\cT{\mathcal{T}}

\def\be{\boldsymbol{\varepsilon}}
\def\bbeta{\boldsymbol{\beta}} 
\def\blambda{\boldsymbol{\lambda}}

\def\bmu{\boldsymbol{\mu}}
\def\bSigma{\boldsymbol{\Sigma}}
\def\bPsi{\boldsymbol{\Psi}}

\title{Asymptotic results with estimating equations for time-evolving clustered data}
\author[1]{Laura Dumitrescu \corref{cor1}}
\address[1]{Victoria University of Wellington, New Zealand}
\cortext[cor1]{Corresponding author}
\ead{Laura.dumitrescu@vuw.ac.nz}

\author[2]{Ioana Schiopu-Kratina}
\address[2]{University of Ottawa, Canada}

\begin{document}
\begin{abstract}
\noindent We study the existence, strong consistency and asymptotic normality  of estimators obtained from estimating functions, that are $p-$dimensional martingale transforms. The problem is motivated by the analysis of evolutionary clustered data, with distributions belonging to the exponential family, and which may also vary in terms of other component series. Within a quasi-likelihood approach, we construct estimating equations, which accommodate different forms of dependency among the components of the response vector and establish multivariate extensions of results on linear and generalized linear models,  with stochastic covariates. Furthermore, we characterize estimating functions which are asymptotically optimal, in that they lead to confidence regions for the regression parameters which are of minimum size, asymptotically. Results from a simulation study and an application to a real dataset are included. 

\end{abstract}

\maketitle

\section{Introduction} 

 Regression models with stochastic covariates are useful in many applications such as the analysis of time series, stochastic recursive approximations, dosage adjustment problems, as well as in learning algorithms and artificial neural networks. A recent application of the latter, which is relevant to our approach, is included in \cite{solares-wei-billings17} and is based on a nonlinear autoregressive exogenous (NARX) model. The NARX model incorporates past values as well as exogenous inputs
$$y_{i} = f(y_{i-1}, y_{i-2}, \ldots, y_{i-q}, u_{i-1}, u_{i-2}, \ldots, u_{i-r}) + \varepsilon_i,$$
where $f$ is a nonlinear mapping, $y_{i}$ and $u_{i}$ are, respectively, the output and the external variable, determined at moment $i.$ Approximating $f({\bf v}_i)$ by $\sum_{l=1}^p \theta_l \phi_l({\bf v}_i),$ with given functions $\phi_l$ and ${\bf v}_i=(y_{i-1}, y_{i-2}, \ldots, y_{i-q}, u_{i-1}, u_{i-2}, \ldots, u_{i-r}),$ for binary data $y_{i},$ in \cite{solares-wei-billings17} the authors propose an algorithm based on the following probability model
$$p_i=\frac{1}{1+\exp[-\sum_{l=1}^p \theta_l \phi_l({\bf v}_i)]}.$$
This method, which combines the logistic function with the NARX representation, is used as a classifier for dynamic binary classification and was shown in simulation studies to perform better than other classification techniques, such as $k-$nearest neighbours or random forest.
We refer to \cite{duflo97} for more applications of stochastic regression, including estimation methods and inference for random iterative models. 

In the present work we consider a general approach for the study of an evolutionary cluster, with new data (having a distribution that belongs to the exponential family) arriving at each time step. Given a collection of past observations, the objective is to ``explain'' their variation in terms of an array of exogenous variables. This approach allows the decomposition of the overall variability of the vector of interest into a time trend and a component that can be explained by exogenous variables and which may give rise to nonstationary time series. In the context of time series, for non-Gaussian data, two alternative modelling approaches, ``{\em observation-driven}'' and ``{\em parameter-driven}'', were discussed in \cite{cox81}. In the latter approach, the autocorrelation is induced via a latent process and was used for modelling time series of counts in e.g. \cite{zeger88} and \cite{davis-wang-dunsmuir99}. A comprehensive review of count time series models is included in \cite{fokianos12}. 

 We consider the {\em observation-driven} approach for multivariate data and study the asymptotic behaviour of a sequence of estimators obtained from an estimating equation based on conditional nonlinear models, for each individual time series.  More precisely, at each moment $i$, we assume that, given the past, the conditional distribution of the $j$th response component, $y_{ij}$, is a function of a linear combination ${\bf x}_{ij}^T\bbeta$. Here, the vector ${\bf x}_{ij}$ may depend on past observations as well as exogenous variables and the objective is to estimate $\bbeta$. 
 
Our approach represents an extension of the stochastic models considered in \cite{lai-wei82} and \cite{chen-hu-ying99}, where an asymptotic theory was established for uni-dimensional linear and generalized linear regression models, respectively. As a function of time, the response is a vector-valued time series with correlated components, whose conditional correlation, given the evolution of time, need not be specified. Our work extends some of the elegant asymptotic results presented in \cite{kaufmann87} to accommodate possible correlations among components of the response vector series. 
 
We remark that, in the related case of longitudinal data, the observations among clusters are assumed to be independent, as they correspond to the situation when several correlated measurements are taken on different subjects. A popular approach is to assume a generalized linear model for the marginals, whereas the correlation among the responses from the same subject is acknowledged, but not modelled. Then, it is well known that estimators obtained from associated estimating equations are consistent regardless of the working correlation structure used and the question of asymptotic efficiency was discussed in \cite{jiang-luan-wang07} and \cite{balan-dumitrescu-schiopu10}. This article shows that a similar property holds when clusters are not independent and are assumed to evolve in time according to a specified pattern, by making inference based on estimating functions which, in this case, are martingale transforms.    
 
The article is organized as follows. We introduce our model assumption in Section \ref{assump-sec} and discuss the special case of independent clusters in Section \ref{longitudinal}. In Section \ref{consist} we give sufficient conditions for strong consistency in Theorem \ref{strong-consistency-th}, 
whose proof is based on a fixed point theorem and a strong law of large numbers for martingales. In Theorem \ref{consistency} we simplify the hypotheses of Theorem \ref{strong-consistency-th} to obtain conditions which are easier to verify.

The study of the limiting distributions is presented in Section \ref{normality}. Theorem  \ref{lemmaCLT} gives a central limit theorem for our estimating functions, which then leads to the main result of the section, Theorem \ref{CLTth}. In the latter we give sufficient conditions for the asymptotic normality of estimators obtained as implicit solutions of the associated estimating equations. 

In Section \ref{optim} we give a characterization of asymptotic optimal sequences of functions (in the sense defined in \cite{heyde97}), which lead to asymptotic confidence region of minimal size. The class of estimating functions that we consider consists of square integrable martingale transforms and we show that not only it includes estimating functions with random coefficients, it also allows for slight misspecification of the regressors (see Corollary \ref{corl}). We expand on this particular type of sensitivity analysis in Section \ref{optim}. 

To evaluate the performance of the proposed method, we present results from a limited simulation study in Section \ref{simulations}. Moreover, we illustrate its applicability on a dataset with average daily wind speed in 2018, taken on several buoys located on the shores of lake Michigan.
   

\section{Model assumptions}
\label{assump-sec}

  Let $y_{ij}$ be a random variable on a probability space, $(\Omega, \mathcal{F}, P_{\boldsymbol{\beta}}),$ $\boldsymbol{\beta} \in \mathcal{T} \subset \mathbb{R}^p$ and $\{\mathcal{F}_i\}_{i \ge 1}$ be an increasing sequence of $\sigma$-fields such that, for each $i \ge 1,$ the $m_i \times 1$ response vector, ${\bf y}_i$  is $\mathcal{F}_i$-measurable and the $p \times 1$ vector  ${\bf x}_{ij}$ is $\mathcal{F}_{i-1}$-measurable, where $j \le m_i,$ $\sup_{i\ge 1} m_i$ and $p$ are bounded. 

 For each component $j,$ at each moment $i,$ the conditional distribution of $y_{ij},$ given  $\mathcal{F}_{i-1}$ is assumed to be a member of the exponential family of distributions or similar. More precisely, we assume that the conditional expectation and variance, given $\mathcal{F}_{i-1}$ are specified in terms of $\bbeta$ through the value $\theta_{ij}={\bf x}_{ij}^T \bbeta:$
$$
	\textrm{E}_{\boldsymbol{\beta}}(y_{ij}| \mathcal{F}_{i-1})  = \mu(\theta_{ij}) = \mu_{ij}(\boldsymbol{\beta}),  \ \
	\textrm{Var}_{\boldsymbol{\beta}}(y_{ij}| \mathcal{F}_{i-1})= \mu'(\theta_{ij}) = \phi \sigma_{ij}^2(\boldsymbol{\beta}),
$$
\noindent where $\mu$ is a differentiable link function, with $\mu'>0$ and $\phi$ is an over-dispersion parameter. Since the estimation of $\boldsymbol{\beta}$ does not involve $\phi,$ for simplicity, we assume that $\phi=1.$ Examples of link functions include the identity, the exponential and the logistic function which are used for continuous responses, count data, and binary observations, respectively. 

 \begin{example}
 {\rm	Consider the case of panel data with binary observations and assume that several non-independent time series are observed simultaneously. For each unit $j$ of a cluster of size $m,$  let $\pi_{ij} = P(y_{ij}=1| y_{i-1 j}, \ldots, y_{1j})$ denote the conditional probability of $y_{ij}=1$ given the past observations $y_{i-1j}, \ldots, y_{1j}.$ 
 	These conditional probabilities are usually not known and, as in the case of autoregressive models for continuous data, they can be estimated using past observations (see \cite{fahrmeir-kaufmann87}, for the case $m=1$). Assuming homogeneity within the cluster, let $\displaystyle{\pi_{ij} = \mu({\bf x}_{ij}^T \bbeta),}$
 	where the vector ${\bf x}_{ij}^T = (1, \  y_{i-1j}, \ \ldots, \ y_{i-l j}, \ z_{ij}^{[1]}, \ldots, \ z_{ij}^{[k]})$ may include $l$ previous observations, as well as $k$ exogenous variables. The assumption leads to a nonhomogenous Markov chain of order $l.$ Other possible formulations could include quadratic or higher order interaction terms. 
 	
 	If the natural link function $\mu$ is chosen as $\mu^{-1}(\pi_{ij}) = {\rm logit}(\pi_{ij}) = \ln(\pi_{ij}) - \ln(1-\pi_{ij}),$ then the binary logit model is obtained, where $\displaystyle{\pi_{ij} = \frac{\exp({\bf x}_{ij}^T \bbeta) }{1+\exp({\bf x}_{ij}^T \bbeta)}}.$
 	Under the additional assumption that the observations within the cluster are independent, at each occasion $i,$ an estimator of $\bbeta$ can be found as a solution of 
 	$$ \sum_{i=1}^n \sum_{j=1}^m {\bf x}_{ij}[y_{ij} - \mu({\bf x}_{ij}^T \bbeta)] = {\bf 0}.$$
}	
 \end{example}

 In the sequel, we do not consider the individual time series to be independent, and since we model separately the conditional distribution for each response component $y_{ij}$, the multivariate conditional distribution of ${\bf y}_i= (y_{i1}, \ldots, y_{im_i})^T,$ given the history $\mathcal{F}_{i-1}$ is not completely specified. Consequently, at each moment in time, a certain form of the conditional correlation structure within the response vector has to be assumed.
	
 Denoting by $\bmu_i(\bbeta) =(\mu({\bf x}_{i1}^T\bbeta), \ldots, \mu({\bf x}_{im_i}^T\bbeta))^T,$ the model assumption imply that the residuals $\be_i={\bf y}_i - {\bmu_i(\bbeta)}$ form a martingale difference sequence with respect to $\mathcal{F}_i,$ that is, $\be_i$ is $\mathcal{F}_i-$measurable and ${\rm E}_{\bbeta}(\be_i|\mathcal{F}_{i-1})={\bf 0},$ for every $i \ge 1.$ It is important to remark that, unlike longitudinal data, where $\be _1, \ldots, \be_i$ are zero mean, independent random vectors, our framework provides a natural generalisation to the situation when there exists a time trend within the cluster ${\bf y}_1, \ldots, {\bf y}_n.$ This also distinguishes our research from the case of transitional models, where the evolution of time is assumed within the components of each of the $n$ independent response vectors.
 
	 Let $\bSigma_{i}(\bbeta)= {\rm Cov}_{\bbeta}({\bf y}_i|\mathcal{F}_{i-1})$ denote the conditional covariance matrix within the cluster, after the time effect is removed and $\bar{{\bf R}}_i(\bbeta)$ be the corresponding conditional correlation matrix. Based on the vector of observations ${\bf y}_1, \ldots, {\bf y}_n$ and ${\bf x}_{11}, \ldots, {\bf x}_{1m_i}, \ldots, {\bf x}_{n 1}, \ldots, {\bf x}_{nm_i},$ if $\bSigma_{i}(\bbeta)$ is specified, an estimator could be computed by solving the score equations ${\bf g}_n^{s}(\boldsymbol{\beta})={\bf 0},$ where
	 \begin{equation}
	 {\bf g}_n^{s}(\boldsymbol{\beta})=\sum_{i=1}^n \left(\frac{\partial \bmu_i(\bbeta)}{\partial \bbeta^T} \right)^T\bSigma_{i}(\bbeta)^{-1} [{\bf y}_i-\bmu_i(\boldsymbol{\beta})]. \label{quasi-score}
	 \end{equation}
  	Let ${\bf A}_i({\bbeta})$ be the diagonal matrix with random entries $\sigma_{i1}^2(\bbeta), \ldots, \sigma_{im_i}^2(\bbeta)$ and so 
		\begin{equation}
		\bar{{\bf R}}_i(\bbeta)={\bf A}_i({\bbeta})^{-1/2}\bSigma_{i}(\bbeta){\bf A}_i({\bbeta})^{-1/2}; \nonumber
		\end{equation}
	in the sequel we characterize properties of estimating equations that use surrogate matrices $\mathcal{R}_{i}(\boldsymbol{\beta})$ as substitutes 
\begin{equation}
	{\bf g}_n(\boldsymbol{\beta})=\sum_{i=1}^n {\bf X}_i^T {\bf A}_i(\boldsymbol{\beta})^{1/2}\mathcal{R}_{i}(\boldsymbol{\beta})^{-1}{\bf A}_i(\boldsymbol{\beta})^{-1/2} [{\bf y}_i-\bmu_i(\boldsymbol{\beta})]. \label{GEE*}
	\end{equation}

The replacement matrices could be selected by the analyst, for example the identity matrix can be chosen as a working correlation matrix within the cluster leading to
\begin{equation}
{\bf g}_n^{\rm ind}(\boldsymbol{\beta}) = \sum_{i=1}^n {\bf X}_i^T ({\bf y}_i - \bmu_i(\boldsymbol{\beta})).  \label{GEE-indep}
\end{equation}
Alternatively, if all $\bar{{\bf R}}_i(\boldsymbol{\beta})$ are the same, the substitution matrix could be an estimator of the true correlation within the cluster, in which cases \eqref{GEE*} defines a {\em pseudo-likelihood},  
with \begin{equation}
\mathcal{R}_n^*(\boldsymbol{\beta})=\frac{1}{n}\sum_{i=1}^n {\bf A}_i(\boldsymbol{ \boldsymbol{\beta}})^{-1/2}({\bf y}_i - \bmu_i( \boldsymbol{\beta}))({\bf y}_i - \bmu_i(\boldsymbol{\beta}))^T {\bf A}_i(\boldsymbol{\beta})^{-1/2}. \nonumber
\end{equation}

 In all cases, our theoretical results confirm that the accuracy to which the proxies approximate the conditional correlation matrices is irrelevant for the existence and strong consistency of the estimators (see the hypotheses of our Theorem \ref{consistency}).  

\section{The special case of independent clusters}
\label{longitudinal}

The proposed methodology for producing estimators of the true parameter $\bbeta_0$ is based on estimating functions and their form is analogous to the estimating equations approach for longitudinal data, i.e. data from independent clusters such as individuals, introduced by \cite{liang-zeger86}. The marginal modelling approach for longitudinal data in \cite{liang-zeger86} specifies the marginal mean and marginal variance collected on occasions $j =1, \ldots, m_i$ as a function of a $p-$dimensional vector of (non-random) covariates, measured at each occasion. The estimator obtained from the generalized estimating equations (GEE) approach is defined as the solution of ${\bf g}_n^{GEE}(\bbeta)={\bf 0},$ where 
\begin{equation}
{\bf g}_n^{GEE}(\bbeta)=\sum_{i=1}^n{\bf X}_i^T {\bf A}_i(\bbeta) \widehat{\bSigma}_i^{-1} ({\bf y}_i - \mu_{i}(\bbeta)) \label{GEE-form}
\end{equation}
and $\widehat{\bSigma}_i$ is an estimator of the $m_i \times m_i$ covariance matrix of ${\bf y}_i.$ 

Under regularity conditions, the GEE estimators are consistent and with ${\bf X}_i$ non-random, we have ${\rm E}_{\bbeta}({\bf g}_n^{GEE}(\bbeta)) = {\bf 0},$ an important property ensuring the consistency of the estimators. The main feature of the GEE approach is that it renders estimators which are consistent even when the correlation (on which $\widehat{\bSigma}_i$ is based) is misspecified and a popular {\em working} correlation matrix is ${\bf I}_{m_i \times m_i},$ i.e. the {\em working independence assumption}.

When covariates are random, some form of conditioning is needed in the specification of the model and it is at this point that our research differs from the GEE approach. In the case of longitudinal models, a direct extension of the marginal modelling approach would be to assume that ${\rm E}_{\bbeta}(y_{ij}| {\bf x}_{ij}) = \mu_{ij}(\bbeta),$ where ${\bf x}_{ij}$ are random vectors which vary in time $j.$ The $k-$th component of a term $i$ in \eqref{GEE-form} is
$$ \sum_{s=1}^{m_i} \sum_{t=1}^{m_i} {x}_{it}^{[k]} a_{it}(\bbeta) \widehat{v}_{ts} (\bbeta) [y_{is} - \mu_{is}(\bbeta)],$$
where $x_{it}^{[k]}$ is the $k-$th component of ${\bf x}_{it},$ $a_{it}(\bbeta)$ is the entry $t$ of the diagonal matrix ${\bf A}_i(\bbeta)$ and $\widehat{v}_{ts}(\bbeta)$ is the $(t,s)$ entry of $\widehat{\bSigma}_i^{-1}.$ While ${\rm E}_{\bbeta}[y_{is} - \mu_{is}(\bbeta)| {\bf x}_{is}] = 0,$ for any $1 \le s \le m_i,$ this conditional expectation cannot traverse factors $x_{it}^{[k]},$ for $t \neq s,$ i.e., in the ``past'' and ``future''. The same remark holds true for {\em transition models} in \cite{pepe-couper97}. In this case, the mean of $y_{ij}$ is conditioned on the entire past history, represented by $y_{is},$ $s < j,$ $x_{is},$ $s \le j.$ We note that models in survival analysis are transition models in this sense.  

It is therefore not guaranteed that the equation of type \eqref{GEE-form} is unbiased for marginal or transition models, with random and time varying covariates and this problem was first brought forth in \cite{pepe-anderson94}. They propose the use of the working independence assumption, or the imposition of a strong condition on the covariates ${\bf x}_i$ (their condition (5)). The authors of \cite{clement-strawderman09} encountered the same problem in the attempt to use the GEE approach in connection with the study of recurrent events and opted for the working independence assumption. It is shown in \cite{lai-small07} that the unrestricted use of the working independence assumption leads to a loss of efficiency in many situations. The authors  divide the covariates into three types and fit estimating equations appropriate to each type. The resulting blocks of estimating functions are then put together by applying the generalized method of moments. 

For longitudinal data, there is often a temporal direction, ascribed to the index $1 \le j \le m_i,$ which is apparent in transition models. In this article, we assume that there is an unspecific undirected correlation among all $m_i$ coordinates of a time series. The evolution in time of the cluster which represents the time series is tracked by the index $i.$ 

\section{Existence and strong consistency}
\label{consist}

For inference, we assume the existence of a true parameter $\bbeta_0$ and hence work on a probability space $(\Omega, \mathcal{F}, P_{\boldsymbol{\beta}_0}),$ which is  $P_{\boldsymbol{\beta}_0}$-complete. The measurability properties of the regressors ensure that our estimating equations are transform martingales, and thus martingale strong convergence results apply. A sufficient condition for strong consistency follows by applying a strong convergence theorem for martingales, as it appears in \cite{lin94}.  
 
  We start by investigating the strong consistency in a linear multivariate framework, with correlated responses. The following example shows that a sufficient condition for consistency requires the minimum eigenvalue of the conditional covariance matrix to converge to infinity faster than the logarithm of the maximum eigenvalue.  More precisely, if $\lambda_{\min}({\bf F}_n)$ tends to infinity faster than $[\log \lambda_{\max}({\bf F}_n)]^{\nu},$ for some $\nu >1,$ then strong consistency of $\widehat{\bbeta}_n$ holds. This result extends Corollary 3 of \cite{lai-wei82} to the multivariate, linear time series context.  
\begin{example}
{\em Consider the multivariate linear regression model with stochastic covariates
$$ {\bf y}_i = {\bf x}_{i}^{[1]} \beta_1 + \ldots {\bf x}_{i}^{[p]} \beta_p + \varepsilon_i, \ i=1, \ldots, m.$$
Here, each of the $m$ rows of ${\bf X}_i = ({\bf x}_{i}^{[1]}, \ldots, {\bf x}_{i}^{[p]})$ is a sequentially determined random vector, that depends on previous values given by the corresponding observations in ${\bf X}_1, {\bf y}_1, \ldots, {\bf X}_{i-1}, {\bf y}_{i-1}.$ Let ${\bf R}_{i}(\alpha)$ denote the assumed correlation structure (possibly depending on a nuisance parameter, $\alpha$) within each cluster ${\bf y}_i,$ and $\bar{{\bf R}}_i$ denote the true conditional correlation. Then, the estimating function ${\bf g}_n(\boldsymbol{\beta})$ is given by 
$${\bf g}_n(\boldsymbol{\beta}) = \sum_{i=1}^n {\bf X}_i^T {\bf R}_{i}(\alpha)^{-1} ({\bf y}_i - {\bf X}_i \boldsymbol{\beta}),$$
and in this case, the explicit solution of the associated estimating equation is $$\widehat{\boldsymbol{\beta}}_n = \left(\sum_{i=1}^n {\bf X}_i^T {\bf R}_{i}(\alpha)^{-1} {\bf X}_i\right)^{-1} \left(\sum_{i=1}^n {\bf X}_i^T {\bf R}_{i}(\alpha)^{-1}{\bf y}_i \right).$$ 

The strong consistency of $\widehat{\boldsymbol{\beta}}_n$ is then equivalent to $$\left(\sum_{i=1}^n {\bf X}_i^T {\bf R}_{i}(\alpha)^{-1}{\bf X}_i\right)^{-1} \left(\sum_{i=1}^n {\bf X}_i^T {\bf R}_{i}(\alpha)^{-1}\boldsymbol{\varepsilon}_i\right) \to 0 \ a.s.,$$
where $\boldsymbol{\varepsilon}_i= {\bf y}_i - {\bf X}_i \boldsymbol{\beta}_0,$ and hence is equivalent to 
$${\bf K}_n^{-1} {\bf g}_n \to 0 \ a.s.,$$
where ${\bf K}_n = \sum_{i=1}^n {\bf X}_i^T {\bf R}_{i}(\alpha)^{-1}{\bf X}_i,$ and ${\bf g}_n= \sum_{i=1}^n {\bf X}_i^T {\bf R}_{i}(\alpha)^{-1}\boldsymbol{\varepsilon}_i \mbox{ is a martingale}.$ 
Applying Theorem 4 of \cite{lin94}, if ${\bf F}_n = \sum_{i=1}^n {\bf X}_i^T {\bf R}_{i}(\alpha)^{-1}\bar{{\bf R}}_i {\bf R}_{i}(\alpha)^{-1}{\bf X}_i$ is almost surely positively definite, then
$${\bf F}_n^{-1} {\bf g}_n \to 0 \ a.s.\mbox{ and in } L^2(P_{\bbeta_0}).$$ The required assumptions are the square integrability of ${\bf g}_n$ and that $[\log \lambda_{\max}({\bf F}_n)]^{\nu} = o(\lambda_{\min}({\bf F}_n)),$ for some $\nu>1.$ 

Note that ${\bf F}_n$ depends on $\bar{{\bf R}}_i$ but if there exists a constant $c>0$ such that $ c \le \min_{i \le n} {\lambda}_{\min} [{\bf R}_{i}(\alpha)],$ then the almost sure convergence holds with the normalizer ${\bf K}_n$ instead.
}
\end{example}

We now consider a general approach and investigate the existence and consistency of a sequence of estimators $\{\widehat{\bbeta}_n\}_{n \ge 1},$ which are derived as solutions of estimating equations, associated with functions of the form ${\bf g}_n(\boldsymbol{\beta})=\sum_{i=1}^n{\bf u}_i(\boldsymbol{\beta}).$ 
Here, ${\bf g}_n(\bbeta)$ is a square integrable and continuously differentiable martingale difference with respect to $\mathcal{F}_n.$ We assume that its quadratic variation matrix, ${\bf V}_n(\bbeta)=\sum_{i=1}^n {\rm E}[{\bf u}_i(\bbeta){\bf u}_i(\bbeta)^T |\cF_{i-1}],$ is almost surely positive definite. This matrix plays an important role in establishing the asymptotic behavior of ${\bf g}_n(\boldsymbol{\beta})$ and it represents the amount of information contained in the past history. Furthermore, if the estimating function ${\bf g}_n(\boldsymbol{\beta})$ is a quasi-score, then ${\bf V}_n(\bbeta)$ gives a generalized form of the Fisher information. 
We denote by
\begin{equation}
\cD_n(\boldsymbol{\beta})=-\frac{\partial {\bf g}_n(\boldsymbol{\beta})}{\partial \boldsymbol{\beta}^T}, \quad  \boldsymbol{\beta} \in \cT. \nonumber
\end{equation}
and when no confusion may arise, we omit writing $\boldsymbol{\beta}_0$ when it appears in the argument of a function, e.g., ${\bf g}_n = {\bf g}_n(\boldsymbol{\beta}_0)$, $n \ge 1$.
The main theorem of this section will follow by an application of an appropriate strong law of large numbers for $p-$dimensional martingales, with random normalizer. We begin with an auxiliary result.
\begin{lemma}
\label{random-norm-martingale-convergence}
If $\lambda_{\min}({\bf V}_n) \to \infty$, then, for any $\delta>0$
\begin{equation}
\frac{{\bf g}_n}{[\lambda_{\max}({\bf V}_n)]^{1/2+\delta}} \longrightarrow 0, \ a.s.\mbox{ and in }L^2(P_{\bbeta_0}). \nonumber
\end{equation} 
\end{lemma}

{\bf Proof.} For any $k \in \{1, \ldots p\}$ we show that  
\begin{equation}
\frac{{\bf g}_n^{[k]}}{[\lambda_{\max}({\bf V}_n)]^{1/2+\delta}} \longrightarrow 0, \ a.s. \mbox{ and in }L^2(P_{\bbeta_0}), \nonumber
\end{equation} 
where ${\bf g}_n^{[k]}$ is the $k$-th component of the vector ${\bf g}_n$. If ${\bf e}_k^T=(0, \ldots, 1, \ldots, 0)$ denotes the $p$-dimensional vector with the $k$-th component equal to 1, then ${\bf g}_n^{[k]}={\bf e}_k^T{\bf g}_n=\sum_{i=1}^n{\bf e}_k^T {\bf u}_i$ and ${\bf v}_n^{[k]}=\sum_{i=1}^n {\rm E}({\bf e}_k^T{\bf u}_i{\bf u}_i^T{\bf e}_k |\cF_{i-1})={\bf e}_k^T{\bf V}_n{\bf e}_k$. Since $\lambda_{\min}({\bf V}_n) \le {\bf v}_n^{[k]}$ 
and $\lambda_{\min}({\bf V}_n) \to \infty$, we have ${\bf v}_n^{[k]} \to \infty$. This together with Theorem 4, (ii) in \cite{lin94} (with $\alpha=1/2+\delta$ and any $\nu>1/(2\delta)$) implies
 \begin{equation}
\frac{{\bf g}_n^{[k]}}{({\bf v}_n^{[k]})^{1/2+\delta}} \longrightarrow 0, \ a.s.\mbox{ and in }L^2(P_{\bbeta_0}). \nonumber
\end{equation}
The conclusion now follows, since ${\bf v}_n^{[k]} \le \lambda_{\max}({\bf V}_n)$. \hfill $\Box$

\vspace{2mm}

In the sequel, we use notation $\|{\bf x}\|$ to denote the Euclidean norm of a vector ${\bf x}$. For a matrix ${\bf A},$ let $\displaystyle{\| {\bf A} \|=\sup_{\|{\bf x}\|=1}\|{\bf A}{\bf x}\|}$ be its spectral norm, $\displaystyle{\|| {\bf A} \||=\sup_{\|{\bf x}\|=1}|{\bf x}^T{\bf A}{\bf x}|}$ be its numerical radius and let $\lambda_{\max}({\bf A})$ (respectively, $\lambda_{\min}({\bf A})$) denote its maximum (respectively minimum) eigenvalue. The following inequality will be useful 
\begin{eqnarray*}
	 \|| {\rm E} ({\bf A}) \|| = \sup_{\|{\bf x}\|=1} |{\bf x}^T {\rm E}({\bf A}) {\bf x}| 
	 \le \sup_{\|{\bf x}\|=1} {\rm E}|{\bf x}^T {\bf A} {\bf x}|  
	 \le \sup_{\|{\bf x}\|=1} {\rm E} \|| {\bf A} \||.
\end{eqnarray*}
In addition, since $\||{\bf A} \|| \le \|{\bf A} \| \le  2 \|| {\bf A}\||$ we have $\| {\rm E} ({\bf A})\| \le 2 {\rm E}\|{\bf A} \|$.

Let $B_r=\{\boldsymbol{\beta} \in \mathcal{T} \ | \  \|\boldsymbol{\beta} - \boldsymbol{\beta}_0\| \le r \}$ and $\partial B_r = \{\bbeta \in \mathcal{T}  \ | \  \|\boldsymbol{\beta} - \boldsymbol{\beta}_0\| = r  \},$ with $r>0$. 

Given a sequence of estimating functions ${\bf g}_n(\bbeta),$ we demonstrate the existence and strong consistency of an estimator $\widehat{\bbeta}_n$ that solves ${\bf g}_n(\bbeta)={\bf 0}$ by showing that for all sufficiently small values of $r$ and for $n$ sufficiently large, the estimating equation has a root in the sphere $B_r.$ The criterion that we use follows from an application of a fixed point theorem (see, e.g. Theorem 6.3.4 of \cite{ortega-rheinboldt70}) and we note that a similar technique was used to prove the strong consistency of estimators (in a more general context) in e.g. Theorem 12.1 of \cite{heyde97} and Theorem 12.1 of \cite{hutton-ogunyemi-nelson91}.

\begin{theorem}
\label{strong-consistency-th}
Assume that the following conditions hold
\begin{eqnarray*}
& (I) & \lambda_{\rm min}({\bf V}_n) \stackrel{a.s.}\longrightarrow \infty, \ n \to \infty. \\
& (S) & \mbox {For any constant } \ \delta>0  \mbox{ and any } \lambda \in \mathbb{R}^p, \|\lambda\|=1, \mbox{ the following  } \\
& & \mbox{ conditions hold with probability 1 } \\
& & (i) \ \  \lim_{r \to 0}\limsup_{n \to \infty}
[\lambda_{\max}({\bf V}_n)]^{-1/2-\delta}\sup_{\boldsymbol{\beta} \in
B_r} \|| \cD_{n}(\boldsymbol{\beta})-\cD_{n} \||=0, \\ 
& & (ii) \ \liminf_{n \to \infty}[\lambda_{\max}({\bf V}_n)]^{-1/2-\delta}\lambda^T\cD_n \lambda \geq c_0, \mbox{ for some constant } c_0>0.
\end{eqnarray*}
\noindent Then, there exists a sequence of random variables $\{ \widehat{\boldsymbol{\beta}}_n \}_{n \ge 1}$ with values in $\cT,$ and a random number $n_0,$ such that

$(a)$ $P({\bf g}_n(\widehat{\boldsymbol{\beta}}_n)=0, \ \mbox{for all} \ n \geq n_0)=1,$

$(b)$ $\widehat{\boldsymbol{\beta}}_n \stackrel{a.s.}\longrightarrow \boldsymbol{\beta}_0,$ as $n \to \infty.$
\end{theorem}

{\bf Proof.} We show that for any (small) value of $r>0,$ almost surely, there exists a random integer $ n_0$ such that for any $n \ge n_0$ 
\begin{equation}
\sup_{\bbeta  \in \partial B_r }(\bbeta - \bbeta_0)^T [\lambda_{\max}({\bf V}_n)]^{-1/2-\delta}{\bf g}_n(\bbeta) <0. \label{fixed-point}
\end{equation}
From here, the conclusion of the theorem can be shown as follows. Let $r>0$ be arbitrarily fixed and consider the events $${\cal E}_n^1(r) = \{ \omega \ | \  \sup_{\bbeta  \in \partial B_r }(\bbeta - \bbeta_0)^T [\lambda_{\max}({\bf V}_n)]^{-1/2-\delta}{\bf g}_n(\bbeta) <0 \}.$$
On each event, ${\cal E}_n^1(r)$ by Brouwer's fixed point theorem, the equation ${\bf g}_n(\bbeta)=0$ has a root in the sphere $\{\boldsymbol{\beta} \in \mathcal{T} \ | \  \|\boldsymbol{\beta} - \boldsymbol{\beta}_0\| \le r \}.$  Then, ${\cal E}_n^1(r) \subseteq {\cal E}_n^2(r),$ where ${\cal E}_n^2(r)$ denotes the event that there exists $\widehat{\bbeta}_n \in B_r,$ with ${\bf g}_n(\widehat{\bbeta}_n)=0$. From \eqref{fixed-point}, $P(\liminf_n {\cal E}_n^1(r))=1$ and so, with probability 1, there exists a random integer $n_0,$ such that ${\bf g}_n(\widehat{\bbeta}_n)=0,$ for $n \ge n_0$.
  
	Part (b) follows from noting that $$P(\limsup_n \{\|\widehat{\bbeta}_n - \bbeta_0\| \ge r\}) \le P(\limsup_n {\cal E}_n^1(r)^C)=0.$$

We now turn to the proof of \eqref{fixed-point}. Let ${\cal E}$ be the event on which $(S)$ holds, with $P_{\bbeta_0}({\cal E})=1.$ Then, for $\omega \in {\cal E},$ and any $\varepsilon>0$ and there exist a random number $r_{\varepsilon}>0$ and a random integer $N_{\varepsilon}\ge 1$ such that for any $r \in (0, r_{\varepsilon}),$ $n \ge N_{\varepsilon},$ and $\lambda \in \mathbb{R}^p,$ $\|\lambda\|=1,$ we have 
$[\lambda_{\max}({\bf V}_n)]^{-1/2-\delta} |\lambda^T [\mathcal{D}_n(\bbeta) - \mathcal{D}_n ]\lambda| < \varepsilon,$ for all $\bbeta \in B_{r},$ and so 
\begin{equation}
[\lambda_{\max}({\bf V}_n)]^{-1/2-\delta} \lambda^T[\mathcal{D}_n(\bbeta) - \mathcal{D}_n] \lambda  >- \varepsilon. \label{deriv-cont}
\end{equation}
 By the Intermediate Value Theorem applied to the function $[\lambda_{\max}({\bf V}_n)]^{-1/2-\delta}{\bf g}_n(\bbeta)$, with $\bbeta \in \partial B_r,$ $r \in (0, r_{\varepsilon})$ and $n \ge N_{\varepsilon},$ there exists $\bbeta_n ^* \in B_r$ such that 
\begin{equation}
[\lambda_{\max}({\bf V}_n)]^{-1/2-\delta}(\bbeta - \bbeta_0)^T ({\bf g}_n(\bbeta) - {\bf g}_n) = -[\lambda_{\max}({\bf V}_n)]^{-1/2-\delta}(\bbeta - \bbeta_0)^T \mathcal{D}_n(\bbeta_n ^*)(\beta - \bbeta_0). \label{IVT}
\end{equation} 
Due to $(S)(ii),$ on ${\cal E},$ there exist a random integer $N',$ and a constant $c_0>0,$ such that, for any $n \ge N'$ we have $[\lambda_{\max}({\bf V}_n)]^{-1/2-\delta}\lambda^T \mathcal{D}_n\lambda \ge c_0.$ This, together with \eqref{deriv-cont}, implies that for any $r \in (0, r_{\varepsilon}),$ there exists $N_{\varepsilon}'=\max\{N_{\varepsilon}, N'\}, $ such that for any $n \ge N_{\varepsilon}'$ and $\bbeta \in \partial B_r$ 
\begin{eqnarray}
[\lambda_{\max}({\bf V}_n)]^{-1/2-\delta}(\bbeta - \bbeta_0)^T \mathcal{D}_n(\bbeta_n ^*)(\bbeta - \bbeta_0) &=& [\lambda_{\max}({\bf V}_n)]^{-1/2-\delta}\lambda^T \mathcal{D}_n\lambda r^2  \nonumber \\
&+& [\lambda_{\max}({\bf V}_n)]^{-1/2-\delta}\lambda^T [\mathcal{D}_n(\bbeta_n^*) - \mathcal{D}_n]\lambda r^2  \nonumber \\
& > & (c_0 - \varepsilon)r^2, \label{eval-term}
\end{eqnarray}
where $\bbeta_n ^* \in B_r$ and we denoted $\lambda=(\bbeta- \bbeta_0)/r.$    

 By Lemma \ref{random-norm-martingale-convergence}, $[\lambda_{\max}({\bf V}_n)]^{-1/2-\delta} {\bf g}_n \to {\bf 0},$ almost surely and let $\bar{\cal E}$ denote the event for which this holds. Then, on $\bar{\cal E},$ for $\varepsilon>0,$ there exists a random integer, $M_{\varepsilon}>0$ such that for any $n \ge M_{\varepsilon},$ we have $[\lambda_{\max}({\bf V}_n)]^{-1/2-\delta} \lambda^T{\bf g}_n < \varepsilon.$ 

Choosing $\varepsilon < c_0/2,$ the latter, together with \eqref{IVT} and \eqref{eval-term} imply that on the event $ {\cal E} \cap \bar{\cal E}$ there exists a random integer $n_0=\max\{N_{\varepsilon}',M_{\varepsilon}\}$ such that for any $n \ge n_0$ 
$$\sup_{\bbeta \in \partial B_r} [\lambda_{\max}({\bf V}_n)]^{-1/2-\delta}(\bbeta - \bbeta_0)^T {\bf g}_n(\bbeta) < (2 \varepsilon - c_0) r^2 <0,$$
which concludes the proof of \eqref{fixed-point}.
\hfill $\Box$

\begin{remark}
{\rm The matrix ${\bf V}_n$ is a form of conditional information which reduces to the Fisher information matrix when clusters ${\bf y}_1, \ldots, {\bf y}_n$ are independent. Furthermore, the assumption $(I)$ has been used in the literature as a sufficient condition for the existence of the maximum likelihood estimator from stochastic processes, see e.g. (6.5) in \cite{hall-heyde80}. The hypothesis $(S)(i)$ requires the sequence $\{[\lambda_{\max}({\bf V}_n)]^{-1/2-\delta}\cD_{n}(\boldsymbol{\beta}) \}_{n \ge 1}$ to be asymptotically equicontinuous at $\bbeta_0,$ with $P_{\bbeta_0}$ equal to 1, for any $\delta>0;$ in \ref{appen} of the Appendix, we provide simplified conditions for this property to hold.
Condition $(S)(ii)$ requires that the normalized derivatives be nonsingular at $\bbeta_0$, so the corresponding normalized estimating functions are locally injective. Overall, condition $(S)$ is required when the regression estimator cannot be expressed explicitly. This is apparent in the proof of Theorem \ref{strong-consistency-th} (see \eqref{eval-term}).
}
\end{remark}

\begin{example}{\rm (Assumptions of Theorem \ref{strong-consistency-th} for the multivariate AR(1) model)} 

	{\rm  Consider an $m-$dimensional first-order autoregressive model, with initial condition ${\bf y}_0={\bf 0}$
		$${\bf y}_i = \beta {\bf y}_{i-1} + \be_i, \ \beta \in \mathcal{T}, $$
		and assume that $\be_i$ are $\mathcal{F}_i-$measurable random variables, with ${\rm E}(\be_i|\mathcal{F}_{i-1})={\bf 0}$ and ${\rm E}(\be_i\be_i^T|\mathcal{F}_{i-1})=\sigma^2 {\bf I}, \ \sigma^2>0,$ where $\mathcal{F}_i=\sigma({\bf y}_1, \ldots, {\bf y}_i).$ The estimating equation $g_n(\beta) = \sum_{i=1}^n {\bf y}_{i-1}^T({\bf y}_i - \beta {\bf y}_{i-1})=0$ gives a solution which corresponds to the classical maximum likelihood estimator of $\beta,$ when $\be_1$ is normally distributed, i.e. $\widehat{\beta}_n = \left(\sum_{i=1}^n {\bf y}_{i-1}^T {\bf y}_{i-1}\right)^{-1} \sum_{i=1}^n {\bf y}_{i-1}^T {\bf y}_i.$ In the case when $m=1,$ it is well known that the estimator is consistent regardless of the distribution of $\be_1.$
		
		Note that ${\rm E}({\bf y}_i|\mathcal{F}_{i-1}) = \beta {\bf y}_{i-1},$  ${\rm E}({\bf y}_i {\bf y}_i^T|\mathcal{F}_{i-1}) = \beta^2 {\bf y}_{i-1}{\bf y}_{i-1}^T+ \sigma^2 {\bf I},$ and so ${\rm Cov}({\bf y}_i| \mathcal{F}_{i-1}) = {\rm Cov}({\be}_i | {\cF}_{i-1}).$ Hence, $\bar{{\bf R}}_{i}={\bf I},$ $i\ge 1,$ ${\bf V}_n = \sigma^2\sum_{i=1}^n {\bf y}_{i-1}^T{\bf y}_{i-1}$ and $\mathcal{D}_n(\beta)=\sum_{i=1}^n {\bf y}_{i-1}^T {\bf y}_{i-1},$ for any $\beta \in \mathcal{T}.$ A sufficient condition for $(I)$ to hold is $\sum_{i=1}^n \|{\bf y}_{i-1}\|^2 \stackrel{a.s.}\longrightarrow \infty.$
		Within conditions $(S)$, $(i)$ is clearly satisfied and under $(I),$ a sufficient condition for assumption $(ii)$ to hold is $\liminf_{n \to \infty}  (\sum_{i=1}^n \|{\bf y}_{i-1}\|^2)^{1/2 - \delta}>0 ,$ for some $\delta>0.$ This condition together with $(I)$ restrict the range of $\delta;$ it has to be no greater than $1/2.$ Hence, the assumption $\sum_{i=1}^n \|{\bf y}_{i-1}\|^2 \stackrel{a.s.}\longrightarrow \infty$ is sufficient for the strong consistency of $\widehat{\beta}_n.$ \hfill $\Box $ 
		} 
\end{example}

The next result is an application of Theorem \ref{strong-consistency-th} to estimating equations of the form \eqref{GEE*}.  Our assumptions are that $\mathcal{R}_{i}(\boldsymbol{\beta})$ are symmetric, positive definite, and have $\mathcal{F}_{i-1}-$ measurable and continuously differentiable entries; these assumptions guarantee that the estimating function ${\bf g}_n(\boldsymbol{\beta})$ is a martingale transform with respect to $\mathcal{F}_n$.	

Assume that $\mu$ is three times continuously differentiable, and for any $r>0$ and $n \geq 1$ let
\begin{eqnarray*}
k_n^{[2]}(r)&=&\sup_{\boldsymbol{\beta} \in B_r}\max_{i \leq n, j \leq
m_i}\left|\frac{\mu''({\bf x}_{ij}^T \boldsymbol{\beta})}{\mu'({\bf x}_{ij}^T
\boldsymbol{\beta})} \right|, \quad k_n^{[3]}(r)=\sup_{\boldsymbol{\beta} \in B_r}\max_{i \leq
n, j \leq m_i}\left|\frac{\mu'''({\bf
x}_{ij}^T \boldsymbol{\beta})}{\mu'({\bf x}_{ij}^T \boldsymbol{\beta})} \right|.
\end{eqnarray*}
As in \cite{xie-yang03} we introduce the assumption 
\begin{eqnarray*}
(AH) \ \mbox{there exist constants} \ C,r_0>0 \mbox{ with}, \  k_{n}^{[l]}(r) \leq C, \ a.s., \ \mbox{for all} \ r \le r_0, \  n \geq 1, \ l=2,3.
\end{eqnarray*}

We denote by ${\bf H}_n'=\sum_{i=1}^n {\bf X}_i^T{\bf A}_i{\bf X}_i.$ The matrix ${\bf H}_n',$ associated with the generalized linear model, corresponds to the design matrix $\sum_{i=1}^n {\bf X}_i^T {\bf X}_i,$ associated with the linear model. Furthermore, ${\bf H}_n'$ replaces ${\bf V}_n$ in Theorem \ref{strong-consistency-th} as appropriate normalizer for the estimating functions \eqref{GEE*} in Theorem \ref{consistency}. 

\begin{theorem}
\label{consistency} Assume that $(AH),$ together with the following conditions, are satisfied
\begin{eqnarray*}
(D) && \ \lambda_{\min}({\bf H}_n') \to \infty \ a.s., \ n \to \infty, \\
(S_{\delta}') && \ \textrm{there exist a random } N \ge 1,  \textrm{ a constant } 0<\delta \le \frac{1}{2} \mbox { and } \mbox{a random } c_0>0 \textrm{ such that} \\
&&\lambda_{\min}({\bf H}_n') \geq c_0 [\lambda_{\max}({\bf H}_n')]^{1/2+\delta} \ a.s.,  \mbox{ for any } n \geq N,\\
(S')(i)&&\lim_{r \to 0} \limsup_{n \to \infty}
[\lambda_{\max}({\bf H}_n')]^{-1/2-\delta}\sup_{\boldsymbol{\beta} \in
B_r}\|| \cD_{n}(\boldsymbol{\beta})-\cD_{n}  \||=0, \ a.s. \\
(H') && \ \mbox{there exists } C>0 \mbox{ such that } \lambda_{\min}(\bar{\bf R}_{n}) > C \ a.s., \mbox{ for any } n\ge 1,\\
(E') && \ \textrm{there exist constants } K, C, \ 0<K \le C \ \textrm{and } K \le \lambda_{\min}(\mathcal{R}_{n}) \le  \lambda_{\max}(\mathcal{R}_{n}) \leq C \ a.s.\\
&& \mbox{for any } n \geq 1.
\end{eqnarray*}
Then there exists a sequence $\{\widehat{\boldsymbol{\beta}}_n\}_{n \ge 1} \subset \mathcal{T}$ and a random number $n_0$ such that 

$(a)$ $P({\bf g}_n(\widehat{\boldsymbol{\beta}}_n)=0, \ \textrm{for all} \ n \geq n_0)=1$;

$(b)$ $\widehat{\boldsymbol{\beta}}_n \stackrel{a.s.} \longrightarrow \boldsymbol{\beta}_0,$ as $n \to \infty.$
\end{theorem}

{\bf Proof.} The proof reduces to the proof of Theorem \ref{strong-consistency-th} given above, once the normalizer is taken to be ${\bf H}_n'$. For instance, in the presence of conditions $(H')$ and $(E')$, conditions $(D)$ and $(I)$ are equivalent. With $0 < \delta <1/2,$ $(S'_{\delta})$ implies condition $(S)(ii)$ of Theorem \ref{strong-consistency-th}. This is so because, for generalized linear models with random regressors, the asymptotically leading term of the normalized derivative ${\cD}_n$ is the normalized matrix ${\bf H}_n'$(see \cite{fahrmeir-kaufmann87} and \cite{kaufmann87}).For non-random regressors, this statement follows from lemmas 4.7--4.9 of \cite{balan-dumitrescu-schiopu10}. \hfill $\Box$

\begin{remark}
{\rm Assumption $(D)$ is a necessary condition for strong convergence (see (1.5) of \cite{lai-wei82}, for the case $m=1,$ $\mu(x)=x$ and the discussion thereafter). Conditions $(S'_{\delta})$ and $(S')(i)$ in the statement of Theorem \ref{consistency} correspond to conditions $(S)(ii)$ and $(S)(i),$ respectively, in Theorem \ref{strong-consistency-th}, with ${\bf H}_n'$ replacing ${\bf V}_n$ as the appropriate normalizer in Theorem \ref{consistency}.
Condition $(H')$ ensures that each matrix $\bar{\bf R}_{n}$ is positively definite and therefore, nonsingular. Condition $(E')$  mimics condition $(H')$ for the approximating matrices $\mathcal{R}_{n}$ which are also positive definite and bounded.  
}
\end{remark}

For all $n \ge 1$, $r>0$, let
\begin{eqnarray*}
\eta_n(r) = \sup_{\boldsymbol{\beta} \in B_r} \max_{i
\leq n, j\le m_i} \left |\left[\frac{\mu'({\bf x}_{ij}^T \boldsymbol{\beta}')}{\mu'({\bf x}_{ij}^T \boldsymbol{\beta})} \right]^{1/2} - 1 \right|. 
\end{eqnarray*}

The proof of Theorem \ref{consistency} gives the following result when the units within the cluster are independent. 
\begin{example}
	
	{\rm Recall that $	{\bf g}_n^{\rm ind}(\boldsymbol{\beta}) = \sum_{i=1}^n {\bf X}_i^T[{\bf y}_i - \bmu_i(\boldsymbol{\beta})] = {\bf 0}$ and assume that $(H'),$ $(D)$ and $(S_{\delta}')$ hold, together with the following condition\\
	 $(K')$ $\lim_{r \to 0} \limsup_{n \to \infty} \eta_n(r) =0. $
		
	There exists then a sequence $\{\widehat{\boldsymbol{\beta}}_n\}_{n \ge 1} \subset \cT$ and a random number $n_0,$ such that
	
	$(a)$ $P({\bf g}^{\rm ind}_n(\widehat{\boldsymbol{\beta}}_n)=0, \ \textrm{for all} \ n \geq n_0)=1;$
	
  $(b)$ $\widehat{\boldsymbol{\beta}}_n \stackrel{a.s.} \longrightarrow \boldsymbol{\beta}_0,$ as $n \to \infty.$
}
\end{example}

\section{Asymptotic normality}
\label{normality}
In what follows, we give sufficient conditions under which a sequence of estimators, obtained as solutions of the estimating equations defined by  \eqref{GEE*}, is asymptotically normally distributed.
The first result gives a Central Limit Theorem (CLT) for martingales ${\bf g}_n = \sum_{i=1}^n {\bf u}_i$ of the form \eqref{GEE*}, with non-stochastic normalizer, 
\begin{equation}
{\bf M}_n^{-1/2} {\bf g}_n  \stackrel{\mathcal{L}}\longrightarrow N({\bf 0}, {\bf I}), \ \mbox{ as } n \to \infty, \nonumber
\end{equation} 
where $\stackrel{\mathcal{L}} \longrightarrow$ denotes the convergence in distribution and 
\begin{equation}
{\bf M}_n = {\rm Cov}( {\bf g}_n) =\sum_{i=1}^n {\rm E}({\bf u}_i{\bf u}_i^{T} ) = \sum_{i=1}^n {\rm E}({\bf X}_i^T {\bf A}_i^{1/2}\mathcal{R}_{i}^{-1} \bar{\bf {R}}_{i} \mathcal{R}_{i}^{-1}{\bf A}_i^{1/2} {\bf X}_i). \nonumber
\end{equation}
We recall that ${\bf V}_n = \sum_{i=1}^n {\rm E}({\bf u}_i{\bf u}_i^{T}|\mathcal{F}_{i-1} ) = \sum_{i=1}^n {\bf X}_i^T {\bf A}_i^{1/2}\mathcal{R}_{i}^{-1} \bar{\bf {R}}_{i} \mathcal{R}_{i}^{-1}{\bf A}_i^{1/2} {\bf X}_i$ and introduce an assumption on its stability: asymptotically, it behaves as the unconditional information matrix\\

{$(A_1)$ \hspace{2mm} ${\bf M}_n^{-1/2} {\bf V}_n^{1/2}  \stackrel{P}\longrightarrow {\bf I}_p $ element-wise, as $n \to \infty.$ \\
	
\begin{remark}
{\rm Assumption $(A_1)$ is referred to as the ergodic case and it has been shown that, in its absence, the limiting distribution may not be normal or may not even exist. For example, when $m_i=1,$ a similar assumption was used in \cite{lai-wei82} to prove the asymptotic normality of the least squares estimator (see their conditions (4.2) and the discussion thereafter).}
\end{remark}

Let $\boldsymbol{\lambda} \in {\mathbb R}^p,$ with $\|\boldsymbol{\lambda}\|=1$ and write
\begin{eqnarray*}
\boldsymbol{\lambda}^T{\bf M}_n^{-1/2} {\bf g}_n = \sum_{i=1}^n \boldsymbol{\lambda}^T {\bf M}_n^{-1/2} {\bf X}_i^T {\bf A}_i^{1/2} \mathcal{R}_{i}^{-1} {\bf A}_i^{-1/2} ({\bf y}_i -\bmu_i)= \sum_{i=1}^n { Z}_{n,i},
\end{eqnarray*}
with ${ Z}_{n,i} = \boldsymbol{\lambda}^T {\bf M}_n^{-1/2} {\bf X}_i^T {\bf A}_i^{1/2} \mathcal{R}_{i}^{-1} {\bf A}_i^{-1/2} ({\bf y}_i -\bmu_i).$ 
Since ${ Z}_{n,i}$ is $\mathcal{F}_{i}-$ measurable, and 
\begin{equation}
{\rm E}({ Z}_{n,i} | \mathcal{F}_{i-1}) := \boldsymbol{\lambda}^T {\bf M}_n^{-1/2}{\bf X}_i^T {\bf A}_i^{1/2} \mathcal{R}_{i}^{-1} {\bf A}_i^{-1/2} {\rm E}({\bf y}_i -\bmu_i | \mathcal{F}_{i-1})  = 0,\nonumber 
\end{equation}
for each $i \le n,$ it follows that $\{{ Z}_{n,i}, \mathcal{F}_{i}\}$ is a martingale difference array.

\begin{theorem}
	\label{lemmaCLT}
	Assume that $(A_1)$, $(E'),$ $(H')$ hold, together with the following
	\begin{eqnarray*}
	&& (I_a) \ \max_{i \le n} \|{\bf M}_n^{-1}({\bf X}_i^T {\bf A}_i {\bf X}_i)\|  \underset{n \to \infty}{ \overset{P}{\longrightarrow}}  0;\\
	&& (N_{\nu}) \ \max_{i\le n}{\rm E}(\|\be_{i}\|^{2+ \nu} | \mathcal{F}_{i-1}) \le C, a.s., \mbox{ for some }\nu>0, \mbox{ where }\be_i = {\bf A}_i^{-1/2} ({\bf y}_i - \boldsymbol{\mu}_i).
	\end{eqnarray*}
	Then, as $n \to \infty,$
	\begin{equation}
	{\bf M}_n^{-1/2} {\bf g}_n  \stackrel{\mathcal{L}}\longrightarrow N({\bf 0}, {\bf I}).
	\end{equation} 
\end{theorem}

{\bf Proof.} Using the Corollary 3.1 in \cite{hall-heyde80}, with $\mathcal{F}_{ni}=\mathcal{F}_i,$ $\eta^2=1,$ $k_n=n$ and the Cramér-Wold theorem, the conclusion follows if we show that
\begin{eqnarray}
\sum_{i=1}^n {\rm E}({ Z}_{n,i}^2 | \mathcal{F}_{i-1}) \stackrel{P}\longrightarrow 1, \mbox{ as } n\to \infty,\label{firststep}
\end{eqnarray}
 
\begin{equation}
\forall \delta>0, \ \sum_{i=1}^n {\rm E}({ Z}_{n,i}^2 1_{\{|Z_{n,i}| \ge \delta\}} | \mathcal{F}_{i-1}) \stackrel{P}  \longrightarrow 0, \mbox{ as } n\to \infty. \label{secondstep}
\end{equation}

Relation \eqref{firststep} is implied by assumption $(A_1)$ since by the definition of $\boldsymbol{\Sigma}_i$ and the $\mathcal{F}_{i-1}-$ measurability property of  ${\bf X}_i,$ ${\bf A}_i,$ $\mathcal{R}_{i},$ we have
\begin{eqnarray*}
	&&\sum_{i=1}^n {\rm E}(\boldsymbol{\lambda}^T {\bf M}_n^{-1/2} {\bf X}_i^T {\bf A}_i^{1/2} \mathcal{R}_{i}^{-1} {\bf A}_i^{-1/2} ({\bf y}_i -\bmu_i)({\bf y}_i -\bmu_i)^T{\bf A}_i^{-1/2} \mathcal{R}_{i}^{-1} {\bf A}_i^{1/2} {\bf X}_i {\bf M}_n^{-1/2} \boldsymbol{\lambda}  |  \mathcal{F}_{i-1}) \\
	&=&\sum_{i=1}^n \boldsymbol{\lambda}^T {\bf M}_n^{-1/2} {\bf X}_i^T {\bf A}_i^{1/2} \mathcal{R}_{i}^{-1} {\bf A}_i^{-1/2} \boldsymbol{\Sigma_i} {\bf A}_i^{-1/2} \mathcal{R}_{i}^{-1} {\bf A}_i^{1/2} {\bf X}_i {\bf M}_n^{-1/2} \boldsymbol{\lambda}   \\
	&=& \lambda^T {\bf M}_n^{-1/2} {\bf V}_n {\bf M}_n^{-1/2}\lambda .
\end{eqnarray*}

We now prove $\eqref{secondstep}$. Let $C>0$ denote a generic constant; using the Cauchy-Schwarz inequality and assumption $(E'),$ almost surely, for any $n \ge 1$ and $i \le n,$ we have the following
\begin{eqnarray}
	{ Z}_{n, i}^2 &=& [\boldsymbol{\lambda}^T {\bf M}_n^{-1/2} {\bf X}_i^T {\bf A}_i^{1/2} \mathcal{R}_{i}^{-1} {\bf A}_i^{-1/2} ({\bf y}_i -\bmu_i)]^2 \nonumber\\
	 &\le& [\boldsymbol{\lambda}^T {\bf M}_n^{-1/2} {\bf X}_i^T {\bf A}_i^{1/2} \mathcal{R}_{i}^{-1} \mathcal{R}_{i}^{-1} {\bf A}_i^{1/2} {\bf X}_i  {\bf M}_n^{-1/2}\boldsymbol{\lambda} ] [\be_i^T \be_i] \nonumber\\
	 &\le&  \lambda_{\max}^2[\mathcal{R}_{i}^{-1}] \|{\bf M}_n^{-1}({\bf X}_i^T {\bf A}_i {\bf X}_i)\| \|\be_i\|^2 \nonumber\\
	 &\le& C \|{\bf M}_n^{-1}  ({\bf X}_i^T {\bf A}_i {\bf X}_i)\| \|\be_i\|^2. \label{ineq1}
\end{eqnarray}

Let $\delta>0.$ 
Then for every $i \le n,$
\begin{equation}
1_{\{|Z_{n,i}| \ge \delta\}} \le 1_{\{ \|\be_i\| \ge C\delta \gamma_{n,i}^{-1}\}}, \label{bound_indicator}
\end{equation} 
where $\displaystyle{\gamma_{n,i}^2:=  \|[{\bf M}_n^{-1} ({\bf X}_i^T {\bf A}_i {\bf X}_i)\|}$ is $\mathcal{F}_{i-1}$ measurable. Condition $(I_a)$ states that $\max_{i\le n}\gamma_{n,i}^2  \underset{n \to \infty}{ \overset{P}{\longrightarrow}}  0.$

As in the proof of \eqref{ineq1}, we obtain, for any $n \ge 1$ and $i \le n,$ a.s.,
\begin{eqnarray*}
{ Z}_{n, i}^2 & \le & 
  \lambda_{\max}[(\bar{{\bf R}}_{i})^{-1}]  \boldsymbol{\lambda}^T {\bf M}_n^{-1/2} {\bf X}_i^T {\bf A}_i^{1/2} \mathcal{R}_{i}^{-1} \bar{{\bf R}}_{i} \mathcal{R}_{i}^{-1} {\bf A}_i^{1/2} {\bf X}_i^T {\bf M}_n^{-1/2} \boldsymbol{\lambda} \|\be_i\|^2 \\
&\le&  C \boldsymbol{\lambda}^T {\bf M}_n^{-1/2} {\bf X}_i^T {\bf A}_i^{1/2} \mathcal{R}_{i}^{-1} \bar{{\bf R}}_{i} \mathcal{R}_{i}^{-1} {\bf A}_i^{1/2} {\bf X}_i^T{\bf M}_n^{-1/2} \boldsymbol{\lambda} \|\be_i\|^2, 
\end{eqnarray*}
where we used $(H')$ for the last inequality. Combining this with \eqref{bound_indicator} leads to
\begin{eqnarray*}
&& \sum_{i=1}^n{\rm E}({ Z}_{n, i}^2 1_{\{|Z_{n,i}| \ge \delta\}}| \mathcal{F}_{i-1})  \\
&&\le \sum_{i=1}^{n} C \boldsymbol{\lambda}^T{\bf M}_n^{-1/2} {\bf X}_i^T {\bf A}_i^{1/2} \mathcal{R}_{i}^{-1} \bar{{\bf R}}_{i} \mathcal{R}_{i}^{-1} {\bf A}_i^{1/2} {\bf X}_i^T {\bf M}_n^{-1/2} \boldsymbol{\lambda} {\rm E}( \|\be_i\|^2 1_{\{ \|\be_i\| \ge C\delta \gamma_{n,i}^{-1}\}}| \mathcal{F}_{i-1}), \\
\end{eqnarray*}
since ${\bf X}_i,$ $\mathcal{R}_{i},$ $\bar{{\bf R}}_{i}$ are 
$\mathcal{F}_{i-1}-$measurable.\\

If for every $i \le n,$ $\|\be_i\| \ge C \delta \gamma_{n,i}^{-1},$ then, for $\nu>0,$ $\|\be_i\|^{\nu} \ge C {\delta}^{\nu} {\gamma_{n,i}}^{-\nu},$ and
\begin{eqnarray}
	C {\delta}^{\nu}\gamma_{n,i}^{-\nu } 1_{\{ \|\be_i\| \ge C \delta \gamma_{n,i}^{-1}\}} \le \|\be_i\|^{\nu}. \nonumber 
\end{eqnarray}
This implies that $C {\delta}^{\nu}\gamma_{n,i}^{-\nu } \|\be_i\|^2 1_{\{ \|\be_i\| \ge C \delta \gamma_{n,I}^{-1}\}} \le \|\be_i\|^{2+{\nu}}$ and, since $\gamma_{n,i}$ is $\mathcal{F}_{i-1}-$measurable, we have
$$ {\rm E}(\|\be_i\|^2 1_{\{ \|\be_i\| \ge C \delta \gamma_{n,i}^{-1}\}} | \mathcal{F}_{i-1}) \le C \delta^{-\nu}\gamma_{n,i}^{\nu} {\rm E}( \|\be_i\|^{2+\nu} | \mathcal{F}_{i-1}), \ i \le n.$$
By assumption $(N_{\nu}),$ the right hand side of the inequality above is almost surely bounded by $ C \delta^{-\nu}\gamma_{n,i}^{\nu} .$

Hence, 
\begin{eqnarray*}
	&& \sum_{i=1}^n{\rm E}({ Z}_{n, i}^2 1_{\{|Z_{n,i}| \ge \delta\}}| \mathcal{F}_{i-1})\\
	&&\le  C \delta^{-\nu} \max_{i\le n}\gamma_{n,i}^{\nu}  \boldsymbol{\lambda}^T{\bf M}_n^{-1/2}\sum_{i=1}^{n} {\bf X}_i^T {\bf A}_i^{1/2} \mathcal{R}_{i}^{-1} \bar{{\bf R}}_{i} \mathcal{R}_{i}^{-1} {\bf A}_i^{1/2} {\bf X}_i^T {\bf M}_n^{-1/2} \boldsymbol{\lambda} \\
	&&=C \delta^{-\nu} \max_{i\le n}\gamma_{n,i}^{\nu} \boldsymbol{\lambda}^T{\bf M}_n^{-1/2} {\bf V}_n {\bf M}_n^{-1/2} \boldsymbol{\lambda}\\
	&&\le C \delta^{-\nu} \max_{i\le n}\gamma_{n,i}^{\nu}\|{\bf M}_n^{-1/2} {\bf V}_n {\bf M}_n^{-1/2} \|.  
\end{eqnarray*}
The expression \eqref{secondstep} follows by $(A_1)$ and $(I_a).$ \hfill $\Box$

\subsection{Proof of Theorem \ref{CLTth}}

Using the Mean Value theorem for ${\bf g}_n,$ for any $n \ge 1,$ on the set $\{{\bf g}_n(\widehat{\boldsymbol{\beta}}_n)=0, \widehat{\boldsymbol{\beta}}_n \in \mathcal{B}_n\}$ we have 
$${\bf g}_n(\widehat{\boldsymbol{\beta}}_n) ={\bf g}_n(\boldsymbol{\beta}_0) - \mathcal{D}_n(\bar{\boldsymbol{\beta}}_n) (\widehat{\boldsymbol{\beta}}_n - \boldsymbol{\beta}_0),$$
where $\bar{\boldsymbol{\beta}}_n \in \mathcal{B}_n$ and so ${\bf M}_n^{-1/2}{\bf g}_n={\bf M}_n^{-1/2} \mathcal{D}_n(\bar{\boldsymbol{\beta}}_n)  (\widehat{\boldsymbol{\beta}}_n - \boldsymbol{\beta}_0).$

Assumption $(A_2)$ implies that $\|{\bf H}_n^{1/2} (\widehat{\boldsymbol{\beta}}_n - \boldsymbol{\beta}_0)\| =O_P(1).$ We write
\begin{eqnarray}
	{\bf M}_n^{-1/2}{\bf g}_n &=& {\bf M}_n^{-1/2} {\bf H}_n^{1/2}[{\bf H}_n^{-1/2} \mathcal{D}_n(\bar{\boldsymbol{\beta}}_n) {\bf H}_n^{-1/2} - {\bf I}]{\bf H}_n^{1/2} (\widehat{\boldsymbol{\beta}}_n - \boldsymbol{\beta}_0) \nonumber\\
	&+&  {\bf M}_n^{-1/2}  {\bf H}_n(\widehat{\boldsymbol{\beta}}_n - \boldsymbol{\beta}_0). \label{CLTdecom}
\end{eqnarray} 
By $(A_3),$ the first term in \ref{CLTdecom} is $o_P(1),$ if $\|{\bf M}_n^{-1/2} {\bf H}_n^{1/2}\|$ is bounded. To prove the last assertion we use the inequality
\begin{eqnarray*}
\boldsymbol{\lambda}^T {\bf V}_n \boldsymbol{\lambda} &\ge& \min_{i \le n} \lambda_{\min}(\bar{{\bf R}}_{i}) \min_{i \le n}    \lambda_{\min}(\mathcal{R}_{i}^{-1}) \boldsymbol{\lambda}^T\sum_{i=1}^n {\bf X}_i^T {{\bf A}_i}^{1/2} \mathcal{R}_{i}^{-1}{{\bf A}_i}^{1/2} {\bf X}_i \boldsymbol{\lambda}, \mbox{ a.s. }
\end{eqnarray*}
Using $(E'),$ $(H')$ and taking expectations gives $ {\bf M}_n^{-1}{\bf H}_n \le C^{-1} {\bf I}_p,$ for some constant $C>0.$
 
 We conclude that 
${\bf M}_n^{-1/2}{\bf g}_n$ and ${\bf M}_n^{-1/2}{\bf H}_n(\widehat{\boldsymbol{\beta}}_n - \boldsymbol{\beta}_0)$ have the same asymptotic distribution. 
An appeal of Theorem \ref{lemmaCLT} concludes the proof. \hfill $\Box$

\begin{remark}
	{\em Condition $(N_{\nu})$ on the conditional moments of the residuals is the same as condition $(N_{\delta})$ of \cite{xie-yang03}. Since ${\bf M}_n = {\rm E}({\bf V}_n),$ then  $(I_a)$ is similar to condition (8) of \cite{xie-yang03} for  bounded cluster sizes, extended here to the case of stochastic covariates.
	 }
\end{remark}
\begin{remark}
	{\em In the case of only one measurement per cluster, i.e. $m_i=1,$ $i\le n$, $(N_{\nu})$ is a Lyapounov type condition. Furthermore, $(I_a)$ was used in \cite{lai-wei82} (see their assumption (4.3)) to obtain the asymptotic distribution of the least squares estimator in stochastic regression. It is analogous to the infinitesimal array condition in Feller's theorem (see also Remark 7 in \cite{xie-yang03} and \cite{billingsley95}) for nonrandom covariates. }
\end{remark}

	For every $n \ge 1,$ we denote by ${\bf H}_n=\sum_{i=1}^n {\rm E}({\bf X}_i^T {{\bf A}_i}^{1/2} \mathcal{R}_{i}^{-1}{{\bf A}_i}^{1/2} {\bf X}_i),$ 
 and $$\mathcal{B}_n = \{\boldsymbol{\beta} \in \mathcal{T}\ | \ \| {\bf H}_n^{1/2}(\boldsymbol{\beta} - \boldsymbol{\beta}_0)\| \le r \}$$
where $r$ is a deterministic constant.

In Theorem \ref{strong-consistency-th}, we gave sufficient conditions for strong  consistency, but here we only assume that there exists a weakly consistent sequence $\{\widehat{\boldsymbol{\beta}}_n\}_{n \ge 1},$  which is a solution of the corresponding estimating equation. Such result can be proved independently of Theorem \ref{strong-consistency-th}, by using a similar technique. 

More precisely, suppose that \\

$(A_2) \ $ there exists a sequence of random variables $\{\widehat{\boldsymbol{\beta}}_n\}_{n \ge 1}$ such that as $n \to \infty, \ $
$P({\bf g}_n(\widehat{\boldsymbol{\beta}}_n) = {\bf 0}) \to 1$ and $P(\widehat{\boldsymbol{\beta}}_n \in \mathcal{B}_n) \to 1.$\\

\begin{theorem}
	\label{CLTth}
	Assume that $(A_1),$ $(A_2),$  $(E'),$ $(H'),$ $(I_a)$ and $(N_{\nu})$ hold,  together with  
	 $$(A_3) \ \sup_{{\boldsymbol{\beta}} \in \mathcal{B}_n}\|{\bf H}_n^{-1/2} \mathcal{D}_n({\boldsymbol{\beta}}) {\bf H}_n^{-1/2} - {\bf I}\| \stackrel{P} \longrightarrow 0, \mbox{ as } n\to \infty.$$
	 
	 Then 
	 ${\bf M}_n^{-1/2}{\bf g}_n$ and ${\bf M}_n^{-1/2}{\bf H}_n(\widehat{\boldsymbol{\beta}}_n - \boldsymbol{\beta}_0)$ have the same asymptotic distribution, and so $${\bf M}_n^{-1/2}{\bf H}_n(\widehat{\boldsymbol{\beta}}_n - \boldsymbol{\beta}_0) \stackrel {\mathcal{L}} \longrightarrow N({\bf 0}, {\bf I}) , \mbox{ as } n \to \infty.$$
\end{theorem}
\begin{remark}
\label{asy-opt-criterion}
{\rm		Under the assumptions of Theorem \ref{CLTth}, for a fixed, $\boldsymbol{\beta}_0 \in \mathcal{T},$ as $n \to \infty,$
		$${\bf M}_n^{-1/2}{\bf H}_n^{}(\widehat{\boldsymbol{\beta}}_n - \boldsymbol{\beta}_0) \stackrel {\mathcal{L}} \longrightarrow N({\bf 0}, {\bf I}).$$
Therefore, to obtain a minimal size asymptotic confidence region, the matrix ${\bf M}_n^{-1/2}{\bf H}_n^{},$ or ${\bf H}_n^{} {\bf M}_n^{-1}{\bf H}_n^{}$ should be maximal with respect to the order of positive semidefinite matrices. Using the concept of optimal estimating equations as presented in \cite{heyde97}, in the next section we characterize asymptotically efficient estimators for models with stochastic covariates.
 }
	\end{remark}
	
\section{Optimal Estimating Equations}
\label{optim}

The estimating function defined in \eqref{quasi-score} 
	is a quasi-score within ${\cal H}_n -$ the family of martingale transforms appearing in \eqref{GEE*}. However, when the intra-cluster correlation is unknown, the use of a proxy for the conditional correlation, $\bar{{\bf R}}_i(\boldsymbol{\beta}),$ renders the associated estimators less efficient. In this section we describe estimating functions of the form \eqref{GEE*}, 
	$${\bf g}_n^*(\boldsymbol{\beta}) = \sum_{i=1}^n {\bf X}_i^T {\bf A}_i(\boldsymbol{\beta})^{1/2}\mathcal{R}_{i}^*(\boldsymbol{\beta})^{-1}{\bf A}_i(\boldsymbol{\beta})^{-1/2} [{\bf y}_i-\bmu_i(\boldsymbol{\beta})]$$ which share, asymptotically, the optimal properties of a quasi-score function, i.e. are asymptotic quasi-score functions and hence the sequences are asymptotically first-order efficient for $\boldsymbol{\beta}$ (as defined in \cite{rao73}) within $\mathcal{H}_n.$ The results are an extension of Theorem 3.9 in \cite{balan-dumitrescu-schiopu10} to the case of stochastic regressors.

We consider here the more general case, when the random covariates are measured with error, on which we do not impose any model assumptions. We require, however, that the errors tend to become negligible as $n$ increases.

Let $\{\boldsymbol{\delta}_i\}_{i \ge 1}$ be a sequence of $p \times m_i$ random matrices, $\|{\boldsymbol{\delta}}_i\| \le d,$ for some $d>0.$ We define ${\bf Y}_i(\boldsymbol{\beta}, \boldsymbol{\delta}_i)^T: = ({\bf X}_i + \boldsymbol{\delta}_i)^T {\bf A}_i(\boldsymbol{\beta}, \boldsymbol{\delta}_i)^{1/2},$ where ${\bf A}_i(\boldsymbol{\beta}, \boldsymbol{\delta}_i)$ is obtained from ${\bf A}_i(\boldsymbol{\beta})$ by substituting ${\bf x}_{ij} + \boldsymbol{\delta}_{ij}$ for ${\bf x}_{ij}$ in $\mu'({\bf x}_{ij}^T\boldsymbol{\beta}),$ $j=1, \ldots, m_i,$ and let ${\bf Y}_i(\boldsymbol{\beta}):={\bf Y}_i(\boldsymbol{\beta}, {\bf 0}),$ $i\ge 1.$ Similarly, when $\mathcal{R}_i^*(\boldsymbol{\beta})$ depends on ${\bf X}_j,$ $j=1, \ldots, i$, we replace ${\bf X}_j$ by ${\bf X}_j + \boldsymbol{\delta}_j$ to obtain $\mathcal{R}_i^*(\boldsymbol{\beta}, \boldsymbol{\delta}_i),$ $i \ge 1.$ We define
\begin{equation}
{\bf g}_n^*(\boldsymbol{\beta}, \boldsymbol{\delta}) := \sum_{i=1}^n {\bf Y}_i(\boldsymbol{\beta}, \boldsymbol{\delta}_i)^T \mathcal{R}_{i}^*(\boldsymbol{\beta}, \boldsymbol{\delta}_i)^{-1} {\bf A}_i(\boldsymbol{\beta}, \boldsymbol{\delta}_i)^{-1/2} ({\bf y}_i - \bmu_i(\boldsymbol{\beta})). \nonumber 
\end{equation}

The next result gives necessary conditions for and shows that, under some continuity and uniform integrability conditions, the optimality property is preserved even if the covariates are recorded with error. We introduce the notation
\begin{eqnarray*}
{\bf H}_n^{\rm ind}(\boldsymbol{\beta}) &=& \sum_{i=1}^n {\rm E}_{\boldsymbol{\beta}}\left[ {\bf X}_i^T {{\bf A}_i}(\boldsymbol{\beta}) {\bf X}_i\right], \\
\bar{{\bf L}}_i(\boldsymbol{\beta}, \boldsymbol{\delta}_i) &=& {\bf Y}_i(\boldsymbol{\beta}, \boldsymbol{\delta}_i)^T \bar{{\bf R}}_i(\boldsymbol{\beta}, \boldsymbol{\delta}_i)^{-1}{\bf Y}_i(\boldsymbol{\beta}, \boldsymbol{\delta}_i), \ 1 \le i \le n,\\
{\bf K}_i^*(\boldsymbol{\beta}, \boldsymbol{\delta}_i) &=& {\bf Y}_i(\boldsymbol{\beta}, \boldsymbol{\delta}_i)^T \mathcal{R}_{i}^*(\boldsymbol{\beta}, \boldsymbol{\delta}_i)^{-1}{\bf Y}_i(\boldsymbol{\beta}, \boldsymbol{\delta}_i), \ 1 \le i \le n.
\end{eqnarray*}

\begin{theorem}
\label{theo_opt}
Let $\{\mathcal{R}_n^*(\boldsymbol{\beta})\}_{n \ge 1}$ be a sequence of $m_n \times m_n$ symmetric, positive definite matrices, and which have $\mathcal{F}_{n-1}-$ measurable and continuously differentiable entries. Assume that $(H')$ holds and that, for all $\boldsymbol{\beta} \in \mathcal{T}$

$(D^*) \quad \lambda_{\min}[{\bf H}_{n}^{\rm ind}(\boldsymbol{\beta})] \to \infty,$

$(R)$  there exists $K(\boldsymbol{\beta})>0$ such that $\inf \lambda_{\min}[\mathcal{R}_n^*(\boldsymbol{\beta})] \ge K(\boldsymbol{\beta}),$ $P_{\boldsymbol{\beta}}$ a.s.,

$(O_1) \quad {\cal R}_{n}^*(\boldsymbol{\beta})- \bar{{\bf R}}_n(\boldsymbol{\beta}) \stackrel{P_{\boldsymbol{\beta}} }{\longrightarrow} 0 \quad \textrm{element-wise},$ as $n \to \infty.$

\indent  In addition, assume that there exists a sequence of $p \times m_i$ random matrices $\{\boldsymbol{\delta}_i\}_{i \ge 1},$ such that

$(O_2) \ \max\{\left\|\mathcal{R}_{i}^*(\boldsymbol{\beta}, \boldsymbol{\delta}_i)^{-1}- \mathcal{R}_{i}^*(\boldsymbol{\beta})^{-1}\right\|, \|{\bf Y}_i(\boldsymbol{\beta}, \boldsymbol{\delta}_i) - {\bf Y}_i(\boldsymbol{\beta})\| \} \le \frac{1}{2^i},$ $i \ge 1$
 
$(O_3) \ \{\|{\bf Y}_i(\boldsymbol{\beta}, \boldsymbol{\delta}_i)\|^2\}_{i \ge 1} \mbox { is } {\rm E}_{\boldsymbol{\beta}} \mbox{- uniformly integrable },$ 

$(O_4) \ \bar{{\bf L}}_i(\boldsymbol{\beta}, \boldsymbol{\delta}_i)^{-1} \mbox{ and } {\bf K}_i^*(\boldsymbol{\beta}, \boldsymbol{\delta}_i)^{-1}, \ i \ge 1 \mbox{ exist and}$  
$$ \inf_{i \ge 1}{\rm E}_{\boldsymbol{\beta}} \{ \lambda_{\min}[ {\bf Y}_i(\boldsymbol{\beta}, \boldsymbol{\delta}_i)^T {\bf Y}_i(\boldsymbol{\beta}, \boldsymbol{\delta}_i)]\} > 0.$$

\noindent Then the sequences $\{{\bf g}_n^*(\boldsymbol{\beta})\}_{n \ge 1}$ and $\{{\bf g}_n^*(\boldsymbol{\beta}, \boldsymbol{\delta})\}_{n \ge 1}$ are asymptotically optimal sequences in $\{\mathcal{H}_n\}_{n \ge 1}$.
\end{theorem}

\begin{remark}
{\rm In condition $(D^*)$ of Theorem \ref{theo_opt}, ${\bf H}_{n}^{\rm ind}$ is the expected value of, and replaces the information matrix ${\bf H}_n'$ in condition $(D)$ of Theorem \ref{consistency}. The use of the expected value in this context is required when the covariates are random. This last remark also motivates the definitions ${\bf M}_n^*(\bbeta)$ and $\bar{\bf M}_{n}(\bbeta)$ (after formula \eqref{Hn^steluta} in the Appendix). Condition $(O_1)$ is condition $(C)$ in \cite{balan-dumitrescu-schiopu10}. Theorem \ref{theo_opt} shows that, when $\bar{{\bf L}}_i(\boldsymbol{\beta}):=\bar{{\bf L}}_i(\boldsymbol{\beta}, {\bf 0})$ or ${\bf K}^*_{i}(\boldsymbol{\beta}):= {\bf K}^*_{i}(\boldsymbol{\beta}, {\bf 0})$ are singular for some $i$ and $\bbeta$ in Theorem 3.9 of \cite{balan-dumitrescu-schiopu10}, one can slightly modify the covariates to obtain non-singular matrices for which the results of Theorem \ref{theo_opt} hold.}
\end{remark}

\begin{corollary}
\label{corl}
Assume that $(H'),$ $(D^*),$ $(R)$ and $(O_1)$ hold. Assume further that ${\mathcal{R}}_i^*(\boldsymbol{\beta})$ is continuous in ${\bf X}_j,$ $1 \le j \le i,$ ${\bf A}_i(\boldsymbol{\beta})$ is continuous in ${\bf X}_i,$ $i \ge 1$ and that $(O_3')$ and $(O_4')$ hold for any $\boldsymbol{\beta} \in \mathcal{T},$

$(O_3') \ \{\|{\bf Y}_i(\boldsymbol{\beta})\|^2\}_{i \ge 1} \mbox { is } {\rm E}_{\boldsymbol{\beta}} \mbox{- uniformly integrable },$ 

$(O_4') \ \inf_{i \ge 1}{\rm E}_{\boldsymbol{\beta}} \{ \lambda_{\min}[ {\bf Y}_i(\boldsymbol{\beta})^T {\bf Y}_i(\boldsymbol{\beta})]\} > 0.$

Then there exists a sequence of random matrices $\{\boldsymbol{\delta}_i\}_{i \ge 1}$ that satisfy the conditions of Theorem \ref{theo_opt} and therefore both $\{{\bf g}_n^*(\boldsymbol{\beta})\}_{n \ge 1}$ and $\{{\bf g}_n^*(\boldsymbol{\beta}, \boldsymbol{\delta})\}_{n \ge 1}$ are asymptotically optimal  sequences in $\{\mathcal{H}_n\}_{n \ge 1}.$
\end{corollary}

{\bf Proof.} The continuity assumptions imply that there exists a sequence of $p \times m_i$ random matrices $\{\boldsymbol{\delta}_i\}_{i \ge 1}$ such that $(O_2)$ holds. These matrices can be chosen such that, using $(O_4'),$ assumption $(O_4)$ is satisfied. By $(O_2),$ the uniform integrability conditions $(O_3)$ and $(O_4')$ are equivalent. Thus, the hypotheses of Corollary \ref{corl} imply the assumptions of Theorem \ref{theo_opt}. \hfill $\Box$

\begin{remark}
{\rm It is always possible to find matrices $\{\boldsymbol{\delta}_i\}_{i \ge 1}$ such that the inverse matrices in $(O_4)$ exist. However, for $\boldsymbol{\delta} \equiv {\bf 0},$ these conditions should be imposed in Theorem 3.9 of \cite{balan-dumitrescu-schiopu10}. The continuity assumption in Corollary \ref{corl} not only ensures that the matrices $(O_4)$ exist, it also allows us to define a large class of sequences that are asymptotic quasi-scores in $\mathcal{H}_n,$ $n \ge 1,$ as long as the norm of $\boldsymbol{\delta}_i$ is ``small'', $i \ge 1.$ This means that even with slightly misspecified regressors ${\bf X}_i,$ $i \ge 1$ in the model, we can still obtain efficient estimators of the regression parameter.}
\end{remark}

The next result gives an example for which hypothesis $(O_1)$ in Theorem \ref{theo_opt} is satisfied. Its proof is an application of a martingale strong law of large numbers, under the assumptions of equal size clusters, constant correlation across clusters and the existence of a higher order moment of normed residuals.  
\begin{prop}
	\label{optimal-ex}
	Assume that $\bar{{\bf R}}_{i}(\boldsymbol{\beta}) :=\bar{{\bf R}}(\boldsymbol{\beta}),$ for any $\boldsymbol{\beta} \in \mathcal{T}$ and $i \ge 1.$ 
	Then, if the following assumption holds\\
		$(N_{\delta}) \ $ ${\rm E}_{\beta} \|{\bf A}_i(\boldsymbol{\beta})^{-1/2} ({\bf y}_i - \boldsymbol{\mu}_i(\boldsymbol{\beta}))\|^{2+\delta_{\beta}} \le C_{\beta},$ for all $i \ge 1$ and some positive constants $\delta_{\beta}$ and $C_{\beta},$ 
		we have 
		\begin{equation}
		\mathcal{R}_{n}^*(\bbeta)- \bar{{\bf R}}(\bbeta) {\longrightarrow}{\bf 0}, \mbox{ a.s. and in } L^1(P_{\bbeta}) \mbox{ element-wise, as } n\to \infty, \label{conv-r-sample}
		\end{equation}
		where
	\begin{equation}
	\mathcal{R}_{n}^*(\boldsymbol{\beta}) = \frac{1}{n} \sum_{i=1}^n {\bf A}_i(\boldsymbol{\beta})^{-1/2} ({\bf y}_i - \boldsymbol{\mu}_i(\boldsymbol{\beta}))  ({\bf y}_i - \boldsymbol{\mu}_i(\boldsymbol{\beta}))^T{\bf A}_i(\boldsymbol{\beta})^{-1/2}. \label{R-optimal}
	\end{equation}
		\end{prop}
		
	{\bf Proof.} Let $\bbeta \in \mathcal{T},$ $1\le j,k \le m$ be arbitrary and denote by $r^{*(n)}_{jk}(\bbeta)$ and $\bar{r}_{jk}(\bbeta)$ the elements of matrices $\mathcal{R}_{n}^*(\boldsymbol{\beta})$ and $\bar{{\bf R}}(\boldsymbol{\beta})$ respectively. Then, we write 
		\begin{eqnarray*}
		r^{*(n)}_{jk}(\bbeta) - \bar{r}_{jk}(\bbeta) &=& \frac{1}{n} \sum_{i=1}^n \sigma_{ij}(\bbeta)^{-1}\left\{ (y_{ij} - \mu_{ij}(\bbeta))(y_{ik} - \mu_{ik}(\bbeta)) \right. \\
		&&\left.   -{\rm E}_{\bbeta}[ (y_{ij} - \mu_{ij}(\bbeta))(y_{ik} - \mu_{ik}(\bbeta)) | \mathcal{F}_{i-1}]\right\} \sigma_{ik}(\bbeta)^{-1}.
		\end{eqnarray*}
With notation $\psi_{jk}^{i}(\bbeta) =  \sigma_{ij}(\bbeta)^{-1}\left\{ (y_{ij} - \mu_{ij}(\bbeta))(y_{ik} - \mu_{ik}(\bbeta)) \right\}\sigma_{ik}(\bbeta)^{-1} \in \mathcal{F}_i,$ Lemma A.1 of \cite{balan-schiopu05} gives that 
$$ \frac{1}{n} \sum_{i=1}^n \left[\psi_{jk}^{i}(\bbeta) - {\rm E}_{\bbeta}[\psi_{jk}^{i}(\bbeta)| \mathcal{F}_{i-1}]\right] \longrightarrow 0, \ a.s. \mbox{ and in }L^{\delta}(P_{\bbeta}),$$
provided that $\sup_{i \ge 1}{\rm E}_{\bbeta}|\psi_{jk}^{i}(\bbeta)|^{\delta} < \infty,$ for some $\delta \in (1, 2].$

The latter is fulfilled due to assumption $(N_{\delta})$ and the fact that $\displaystyle{ \max_{1 \le j,k \le m} |\psi_{jk}^{i}(\bbeta)| \le \| \Psi_i(\bbeta)\|},$ where $\Psi_i(\bbeta) :={\bf A}_i(\boldsymbol{\beta})^{-1/2} ({\bf y}_i - \boldsymbol{\mu}_i(\boldsymbol{\beta}))  ({\bf y}_i - \boldsymbol{\mu}_i(\boldsymbol{\beta}))^T{\bf A}_i(\boldsymbol{\beta})^{-1/2},$ defines the matrix with components $\psi_{jk}^{i}(\bbeta).$ \hfill $\Box$
	
\begin{example}
	\label{R-pseudo-optim}
	{\rm A {\em pseudo-likelihood} estimator can be defined using 
		$$\widetilde{\mathcal{R}}_n = \mathcal{R}_{n}^*(\widehat{\bbeta}_n^{\rm ind}),$$
		where $\mathcal{R}_{n}(\bbeta)$ is the average of squared conditional residuals in  \eqref{R-optimal} and $\widehat{\bbeta}_n^{\rm ind}$ is an initial estimator obtained from solving equation \eqref{GEE-indep}.
		Then, an (asymptotically) efficient estimator of $\bbeta$ can be found by solving the estimating equation
		\begin{equation}
		\widetilde{{\bf g}}_n(\boldsymbol{\beta})=\sum_{i=1}^n {\bf X}_i^T {\bf A}_i(\boldsymbol{\beta})^{1/2}\widetilde{\mathcal{R}}_{i-1}^{-1}{\bf A}_i(\boldsymbol{\beta})^{-1/2} ({\bf y}_i-\bmu_i(\boldsymbol{\beta}))={\bf 0}.
		\end{equation}
		As noted in Example 2.3 of \cite{balan-dumitrescu-schiopu10}, under the assumption that $\bar{{\bf R}}_{i}(\boldsymbol{\beta}) :=\bar{{\bf R}}(\boldsymbol{\beta}),$ for any $\boldsymbol{\beta} \in \mathcal{T}$ and $i \ge 1,$ the analogous of \eqref{conv-r-sample} can be shown for $\widetilde{\mathcal{R}}_n.$ We remark that a similar version of this estimating equation was considered in \cite{balan-schiopu05}, for non-stochastic regressors and, at each step, $i,$ the estimator $\widetilde{\mathcal{R}}_{n}$ was used instead. 
	}

\end{example}

\section{Numerical results}
\label{simulations}
In this section we first examine the convergence and efficiency of the estimators, obtained as solutions of associated martingale transforms estimating equations, with different conditional correlation patterns, and then present an application to real data. 

Consider the simple case of a multivariate second-order autoregressive model specified as 
\begin{equation}
{\bf y}_{i} = {\bf X}_i {\bbeta}_0 + \be_i, \label{multivariateAR1}
\end{equation}
with initial conditions ${\bf y}_{0}$ and ${\bf y}_{1},$ ${\bf X}_i = ({\bf y}_{i-1}, {\bf y}_{i-2})$ and 2-dimensional unknown parameter $\bbeta_0.$ The responses are $m-$dimensional vectors with possibly correlated components, given the time and each ${\bf y}_i$ is assumed to follow a continuous multivariate distribution. Here, the link is the identity function, and let the scale factor, $\phi$ be equal to 1. We consider several choices for the conditional correlation of ${\bf y}_i$, given past observations ${\bf y}_0, {\bf y}_1, \ldots, {\bf y}_{i-1},$ such as 
\begin{enumerate} 
	\item {\em Independence}:  $\bar{{\bf R}}^1={\bf I}_m.$ 
	\item {\em Compound symmetry or exchangeable} (CS): 
$	\bar{{\bf R}}^2= \begin{bmatrix}
1      & \alpha & \ldots & \alpha \\
\alpha & 1      & \ldots & \alpha \\
\vdots & \vdots & \ddots & \vdots \\
\alpha & \alpha & \ldots & 1      \\
\end{bmatrix}.$ 
	\item {\em First-order autoregressive} AR(1):  $\bar{{\bf R}}^3= \begin{bmatrix}
1            & \alpha       & \alpha^2     & \ldots & \alpha^{m-1} \\
\alpha       & 1            & \alpha       & \ldots & \alpha^{m-2} \\
\vdots       & \vdots       & \vdots       & \ddots & \vdots \\
\alpha^{m-1} & \alpha^{m-2} & \alpha^{m-3} & \ldots & 1      \\
\end{bmatrix}.$
\end{enumerate}
The errors are assumed to be independent, and identically distributed $m-$dimensional random vectors, $\be_i\sim N_{m}({\bf 0}, \bSigma),$ with $\bSigma=\bar{{\bf R}}.$ 

To obtain our simulation results, we used the R software (\cite{r-cite20}) as follows. For each $r=1,2, 3, $ we generate values of ${\bf y}_1, \ldots, {\bf y}_n,$ sequentially, according to the model \eqref{multivariateAR1}, with $\be_i\sim N_{m}({\bf 0}, \bar{{\bf R}}^r)$ and then solve the estimating equation
$$\sum_{i=1}^n {\bf X}_{i}^T \widehat {\mathcal{R}}_{i-1}(\bbeta)^{-1}( {\bf y}_i -{\bf X}_i {\bbeta}) ={\bf 0},$$
using several estimators $\widehat {\mathcal{R}}_{i-1}(\bbeta)$ as discussed in what follows.

\vspace{2mm}

Firstly, we use the specified working correlation structures: independence, compound symmetry and AR(1), and denote the resulting solutions of the above equation as  $\widehat{\bbeta}^{(1)}_n,$ $\widehat{\bbeta}^{(2)}_n,$ $\widehat{\bbeta}^{(3)}_n,$ respectively. In each case, we choose $\alpha= 0.7.$ 

\vspace{2mm}

Secondly, we consider the estimator $\widetilde{\bbeta}_n$ defined from a two-step procedure. We let $\widetilde{\mathcal{R}}_{0} = {\bf I}_m,$ $\widetilde{\mathcal{R}}_{1} ={\bf I}_m$ and for each $i > 2$ we proceed as follows. 

$(a)$ Based on the preliminary value,  $\widehat{\bbeta}_{i-1}^{\rm ind} = \left(\sum_{t=1}^{i-1} {\bf X}_{t}^T {\bf X}_{t} \right)^{-1} \sum_{t=1}^{i-1} {\bf X}_{t}^T {\bf y}_t,$ we compute 
	$$\widetilde{\mathcal{R}}_{i-1} = \frac{1}{i-1}\sum_{l=1}^{i-1} ({\bf y}_l - {\bf X}_{l} \widehat{\bbeta}_{i-1}^{\rm ind}) ({\bf y}_l - {\bf X}_{l}\widehat{\bbeta}_{i-1}^{\rm ind})^T.$$

$(b)$ Then, obtain $\widetilde{\bbeta}_n$ by solving the system of estimating equations
	$\displaystyle{\sum_{i=1}^n{\bf X}_{i}^T \widetilde{\mathcal{ R}}_{i-1}^{-1}({\bf y}_i - {\bf X}_{i} \bbeta)={\bf 0}},$ which, in this case, has the closed-form solution
	\begin{equation}
	\widetilde{\bbeta}_n = \left(\sum_{i=1}^n{\bf X}_{i}^T \widetilde{\mathcal{ R}}_{i-1}^{-1}{\bf X}_{i} \right)^{-1}\sum_{i=1}^n{\bf X}_{i}^T \widetilde{\mathcal{ R}}_{i-1}^{-1}{\bf y}_i. \label{beta-pseudo}
	\end{equation}

We generate $s=500$ data sets with $\bbeta_0=(0.5, 0.2)^T,$ $m=5,$ $n=500$ and initial conditions ${\bf y}_0 = {\bf 0}_m,$ ${\bf y}_1 = {\bf 0}_m.$ Tables \ref{biasresults} and \ref{mseresults} compare the finite sample performances of $\widetilde{\bbeta}_n$ with those of $\widehat{\bbeta}^{(j)}_n,$ with $j=1,2,3.$ For each estimator that we consider, we compute estimates for its bias, mean squared error (MSE) and relative efficiency (RE) as follows
\begin{eqnarray*}
	\mbox{Bias}(\widehat{\bbeta}_{n,k}) =\frac{1}{s} \sum_{l=1}^{s} (\widehat{\bbeta}_{n,k}^l - \bbeta_{0,k}), \ \mbox{RB}(\widehat{\bbeta}_{n,k}) &=& \frac{\mbox{Bias}(\widehat{\bbeta}_{n,k})}{\bbeta_{0,k}},
\end{eqnarray*}
\begin{eqnarray*}
	\mbox{MSE}(\widehat{\bbeta}_{n,k}) = \frac{1}{s} \sum_{l=1}^{s} (\widehat{\bbeta}_{n,k}^l - \bbeta_{0,k})^2, \ \mbox{RE}(\widehat{\bbeta}_{n,k}) &=& \frac{\mbox{MSE}(\widehat{\bbeta}_{n,k})}{\mbox{MSE}(\widehat{\bbeta}_{n,k}^{\mbox{true}})}.
\end{eqnarray*}
In the above, the subscript $k \in \{1,2\}$ denotes the $k$th component of the corresponding bivariate vector, whereas the superscript $l \in \{1, \ldots, s\}$ is indicates the estimate obtained at the $l$th simulation run. Moreover, $\widehat{\bbeta}_n^{\mbox{true}}$ is defined as the estimator obtained from the quasi-score \eqref{quasi-score}.

\begin{table} [h!]
\caption{Comparison of simulated relative biases of estimators $\widehat{\bbeta}^{(j)}_n,$ with $j=1,2,3$ and $\widetilde{\bbeta}_n$	\label{biasresults}
 }
\centering
			\begin{tabular}[c]{c l r r r}   
				\hline \hline \\
				Estimator                      &  & $\bar{{\bf R}}^1$ & $\bar{{\bf R}}^2 (\alpha_0=0.7)$ & $\bar{{\bf R}}^3(\alpha_0=0.7)$ \\ 
				\hline
				$\widehat{\bbeta}_{n,1}^{(1)}$ &  & $ 1.3 \cdot 10^{-3}$ & $-2 \cdot 10^{-3}$ & $-0.3 \cdot 10^{-3}$ \\
        $\widehat{\bbeta}_{n,2}^{(1)}$ &  & $-4.8 \cdot 10^{-3}$ & $-11\cdot 10^{-3}$ & $2.2  \cdot 10^{-3}$ \\
	      \hline
        $\widehat{\bbeta}_{n,1}^{(2)}$ &  & $   2 \cdot 10^{-3}$ & $-0.4 \cdot 10^{-3}$ & $-2.8 \cdot 10^{-3}$ \\
        $\widehat{\bbeta}_{n,2}^{(2)}$ &  & $-1.4 \cdot 10^{-3}$ & $-  4 \cdot 10^{-3}$ & $ 8.5 \cdot 10^{-3}$ \\
	      \hline
        $\widehat{\bbeta}_{n,1}^{(3)}$ &  & $1.9 \cdot 10^{-3}$ & $0.8  \cdot 10^{-3}$ & $-2.2 \cdot 10^{-3}$ \\
        $\widehat{\bbeta}_{n,2}^{(3)}$ &  & $1.9 \cdot 10^{-3}$ & $-6   \cdot 10^{-3}$ & $ 9.5 \cdot 10^{-3}$ \\
	      \hline
        $\widetilde{\bbeta}_{n,1}$ &  & $1.9  \cdot 10^{-3}$ & $-0.6 \cdot 10^{-3}$ & $-3.2 \cdot 10^{-3}$ \\
        $\widetilde{\bbeta}_{n,2}$ &  & $-5.5 \cdot 10^{-3}$ & $-2.2 \cdot 10^{-3}$ & $  11 \cdot 10^{-3}$ \\ 
	      \hline
\end{tabular}
\end{table}

\begin{table} [h!]
\caption{Comparison of simulated relative efficiencies of estimators $\widehat{\bbeta}^{(j)}_n,$ with $j=1,2,3$ and $\widetilde{\bbeta}_n$ 	\label{mseresults}}
\centering
			\begin{tabular}[c]{c l r r r}   
				\hline \hline	\\
				Estimator                      &  & $\bar{{\bf R}}^1$ & $\bar{{\bf R}}^2 (\alpha_0=0.7)$ & $\bar{{\bf R}}^3(\alpha_0=0.7)$ \\ 
				\hline
				$\widehat{\bbeta}_{n,1}^{(1)}$ &  & $1$       & $2.87$ & $2.23$ \\
        $\widehat{\bbeta}_{n,2}^{(1)}$ &  & $1$       & $2.90$ & $2.34$ \\
	      \hline
        $\widehat{\bbeta}_{n,1}^{(2)}$ &  & $1.28$    & $1$    & $1.49$ \\
        $\widehat{\bbeta}_{n,2}^{(2)}$ &  & $1.20$    & $1$    & $1.45$ \\
	      \hline
        $\widehat{\bbeta}_{n,1}^{(3)}$ &  & $1.63$    & $1.13$ & $1$ \\
        $\widehat{\bbeta}_{n,2}^{(3)}$ &  & $1.62$    & $1.17$ & $1$ \\
	      \hline
        $\widetilde{\bbeta}_{n,1}$ &  & $1.10$    & $1.08$ & $1.11$ \\
        $\widetilde{\bbeta}_{n,2}$ &  & $1.10$    & $1.16$ & $1.09$ \\ 
	      \hline
\end{tabular}
\end{table}

On the columns of Table \ref{biasresults} we include the results of the relative bias of each component of the estimators that we consider, for different conditional correlation structures that were used to generate the data. We remark that in all scenarios, the bias of each estimator is very small, and this result is consistent with the conclusion of Section \ref{consist}; the estimators are consistent regardless of the conditional correlation pattern chosen in the estimating function. 

Table \ref{mseresults} presents the simulated efficiency of each estimator, relative to that of the estimator obtained as if the conditional correlation were known. In all scenarios, the proposed estimator, $\widetilde{\bbeta}_n$ has smaller MSE, compared to the other estimators. The largest loss of efficiency is seen, as expected, in the case of $\widehat{\bbeta}_n^{(1)}$ when data was generated, using the (CS) or the AR(1) patterns. However, we also note that the results corresponding to $\widehat{\bbeta}_n^{(2)}$ and $\widehat{\bbeta}_n^{(3)}$ are the best possible since, in our computations, we did not estimate the conditional correlation coefficient, $\alpha,$ but rather used its true value in their computation.   

\subsection{Application to the wind speed dataset}
In this section we provide an application to a dataset consisting of daily sensor measurements under different environmental conditions collected by three buoys in lake Michigan in 2018. Two of them are located in Chicago area (\verb!cmti2! and \verb!cnii2!), while the other (\verb!hlnm4!) is on the east side of the lake. The observations were extracted from \verb!data_buoy! available in the R package \verb!forecastML! (\cite{forecastML20}), which, in turn, was obtained from NOAA’s National Buoy Data Center available at \url{https://www.ndbc.noaa.gov/} using the \verb!rnoaa! package. The original dataset contains measurements of fourteen buoys on lake Michigan, recorded from 2012 until 2018 but since many of these have missing values, we chose an appropriate subset containing complete observations of the outcome.   

The variable of interest is considered to be the average daily wind speed in the south-east part of Lake Michigan and a graph of daily recordings from 2018 is included in Figure \ref{graphs}.

\begin{figure}[h!]
\centering
\caption{Daily wind speed recordings in 2018 from buoys located on lake Michigan}
\includegraphics[scale=0.5]{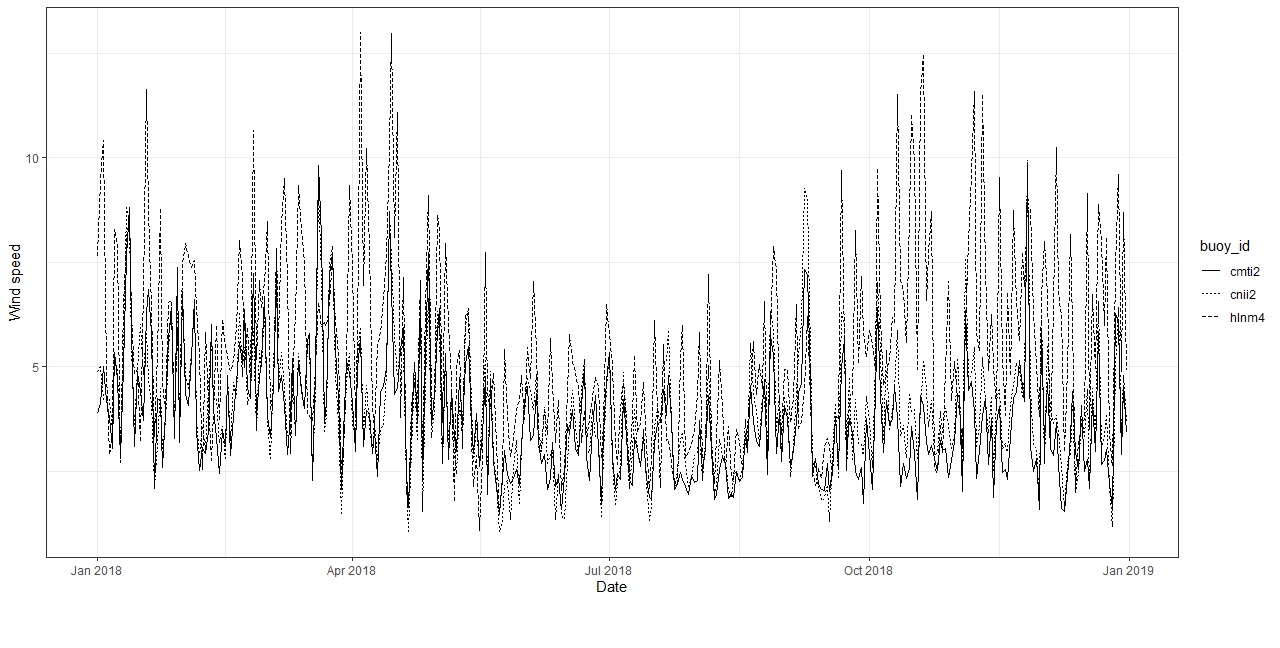}
\label{graphs}
\end{figure} 

As expected, due to their close location, the patterns of variation of wind speed measurements recorded at \verb!cmti2! and \verb!cnii2! are quite similar, whereas measurements from \verb!hlnm4! are higher, on average, and also exhibit a larger variation (we refer to the \verb!data_buoy! dataset for their exact location coordinates).

We use a multivariate model with intercept, random and fixed covariates given by 
$${\bf y}_i= {\bf X}_i\bbeta_0 +\be_i, \ i=0, \ldots, 364. $$
Here ${\bf X}_i=(\mathbbm{1}, {\bf y}_{i-1}, {\bf y}_{i-2}, {\bf z}_{i}),$ with $\mathbbm{1}=(1,1,1)^T,$ $\bbeta_0=(\beta_{0,0}, \beta_{0,1}, \beta_{0,2}, \beta_{0,3})^T$ and where ${\bf y}_{i}=(y_{i1}, y_{i2}, y_{i3})^T$ and ${\bf z}_i=(z_{i1}, z_{i2}, z_{i3})^T$ record the wind speed and air temperature on the $i$th day, at each of the three locations. Hence, $n=364,$ $m=3$ and initial conditions are specified by the recordings on January 1 and January 2 of 2018 (${\bf y}_0$ and ${\bf y}_1$). We note that there were three missing recordings of the air temperature which were replaced by each of their nearest-neighbour values.

Using \eqref{beta-pseudo}, the pseudo-likelihood estimator is  $$\widetilde{\bbeta}=(2.072, 0.441, 0.091, -0.011)^T.$$ 

To obtain a measurement of the accuracy of estimating the components of $\widetilde{\bbeta},$ we apply Theorem \ref{CLTth} and, with notation $\boldsymbol{\Psi} = {\bf H}_n^{-1} {\bf M}_n {\bf H}_n^{-1},$ we write
$${\bPsi}^{-1/2}(\widetilde{\bbeta} - \bbeta_0) \approx N({\bf 0}, {\bf I}).$$
The normalisation terms in $\bPsi$ depend on the unknown correlation matrix $\bar{\bf R}_i,$ but also on the unknown value of $\bbeta_0;$ we approximate $\boldsymbol{\Psi}$ by $\widetilde{\bPsi}= \widetilde{{\bf H}}_n^{-1} \widetilde{{\bf M}}_n \widetilde{{\bf H}}_n^{-1},$ where 
\begin{eqnarray*}
\widetilde{{\bf H}}_n &=& \sum_{i=1}^n {\bf X}_i^T\widetilde{\cal R}_{i-1}^{-1}{\bf X}_i,\\
\widetilde{{\bf M}}_n &=& \sum_{i=1}^n {\bf X}_i^T\widetilde{\cal R}_{i-1}^{-1}({\bf y}_i - {\bf X}_i \widetilde{\bbeta})({\bf y}_i - {\bf X}_i \widetilde{\bbeta})^T\widetilde{\cal R}_{i-1}^{-1}{\bf X}_i.
\end{eqnarray*}  
If $\bPsi$ behaves asymptotically like a constant matrix, then, based on $\widetilde{\bPsi},$ for any $\blambda \in \mathbb{R}^2$ with $\|\blambda\|=1,$ we have  
$$P \left(\blambda^T \bbeta_0 \in (\blambda^T \widetilde{\bbeta}- c_{\alpha/2} \blambda^T \widetilde{\bPsi} \blambda, \blambda^T \widetilde{\bbeta}+ c_{\alpha/2} \blambda^T \widetilde{\bPsi} \blambda) \right) \approx 1-\alpha.$$ 
Here, $c_{\alpha/2}$ denotes the $(1-\alpha)100\%$ upper quantile of the standard normal distribution and the approximation provides, in particular, large sample confidence intervals for the components of $\bbeta_0.$ 

For the wind speed data, the square-roots of the elements on the diagonal of $\widetilde{\bPsi}$ are given by $(0.258, 0.039, 0.038, 0.011).$ Then, the approximate 95\% confidence intervals for $\beta_{0,0},$ $\beta_{0,1},$ $\beta_{0,2},$ and $\beta_{0,3}$ are $(1.566, 2.578),$ $(0.365, 0.517),$ $(0.017, 0.166)$ and $(-0.033, 0.011),$ respectively. 

Since the air temperature is not significant for modelling the wind speed, at $5\%,$ we drop it and fit again the model with ${\bf X}_i=(\mathbbm{1}, {\bf y}_{i-1}, {\bf y}_{i-2})$ to obtain  
$$\widetilde{\bbeta}=(2.386, 0.455, 0.109)^T \mbox{ and }\widetilde{\bPsi}= \left[\begin{array} {rrr} 0.022 & -0.002  & -0.003 \\ -0.002 & 0.001 & -0.001 \\ -0.003 & -0.001 & 0.001
 \end{array} \right].$$ 
These give the following approximate $95\%$ confidence intervals for $\beta_{0,0},$ $\beta_{0,1}$ and $\beta_{0,2}:$ $(2.095, 2.677),$ $(0.393, 0.517)$ and $(0.047, 0.171),$   respectively.

Furthermore, using $\widetilde{\bbeta}$ the average wind speed on January 1, 2019,  in the area where the three buoys are located, can be predicted as 
$$\widehat{\bf y}^{\rm pseudo}= {\bf X}_{365}\widetilde{\bbeta}=(4.444, 4.586, 5.573)^T, \ {\bf X}_{365}=(\mathbbm{1}, {\bf y}_{364}, {\bf y}_{363}).$$ 
Here ${\bf y}_{364}$ and ${\bf y}_{363}$ denote the wind speed average speed on December 31, 2018 and December 30, 2018, respectively.


\newpage
\section{Appendix}
\label{Appendix}
\subsection{Proof of Theorem \ref{theo_opt}}

 Remark \ref{asy-opt-criterion} shows that an asymptotically optimal estimating function can be defined as the one for which ${\bf H}_n^{*}(\bbeta) {\bf M}_n^{*}(\bbeta)^{-1}{\bf H}_n^{*}(\bbeta)$ is maximal within the collection $\{\mathcal{H}_n\}_{n \ge 1}$ for a fixed $\boldsymbol{\beta} \in \mathcal{T}$ and  $n$ large enough, where 
\begin{eqnarray}
{\bf H}_n^*(\boldsymbol{\beta}) &:=& -\textrm{E}_{\boldsymbol{\beta}}\left[\frac{\partial{\bf g}_n^*(\boldsymbol{\beta})}{\partial \boldsymbol{\beta}^T}\right] = \sum_{i=1}^n {\rm E}_{\boldsymbol{\beta}}\left[ {\bf X}_i^T {{\bf A}_i}(\boldsymbol{\beta})^{1/2} \mathcal{R}_{i}^*(\boldsymbol{\beta})^{-1} {{\bf A}_i}(\boldsymbol{\beta})^{1/2} {\bf X}_i\right], \nonumber \\
	                              &:=& \sum_{i=1}^n {\bf K}_i^*(\boldsymbol{\beta}),  \label{Hn^steluta}\\
{\bf M}_n^*(\boldsymbol{\beta}) &:=& \textrm{Cov}_{\boldsymbol{\beta}}\left[{\bf g}_n^*(\boldsymbol{\beta})  \right] \nonumber \\
	                              &=& \sum_{i=1}^n {\rm E}_{\boldsymbol{\beta}}\left[ {\bf X}_i^T {{\bf A}_i}(\boldsymbol{\beta})^{1/2} \mathcal{R}_{i}^*(\boldsymbol{\beta})^{-1} \bar{{\bf R}}_{i}(\boldsymbol{\beta})  \mathcal{R}_{i}^*(\boldsymbol{\beta})^{-1} {{\bf A}_i}(\boldsymbol{\beta})^{1/2} {\bf X}_i\right], \nonumber 
\end{eqnarray}

Furthermore, we employ notations
{\begin{eqnarray*}
\bar{{\bf M}}_n(\boldsymbol{\beta})      &:=& \textrm{Cov}_{\boldsymbol{\beta}}\left[{\bf g}^s_n(\boldsymbol{\beta})  \right] 
                                 =\sum_{i=1}^n {\rm E}_{\boldsymbol{\beta}}\left[{\bf X}_i^T {{\bf A}_i}(\boldsymbol{\beta})^{1/2} \bar{{\bf R}}_{i}(\boldsymbol{\beta})^{-1}	  {{\bf A}_i}(\boldsymbol{\beta})^{1/2} {\bf X}_i \right] \nonumber \\
                                        &:=& \sum_{i=1}^n \bar{{\bf L}}_i(\boldsymbol{\beta}), \nonumber \\  
\bar{{\bf M}}_n(\boldsymbol{\beta}, \boldsymbol{\delta}) &:=& \sum_{i=1}^n \bar{{\bf L}}_i(\boldsymbol{\beta}, \boldsymbol{\delta}_i), \ \bar{{\bf L}}_i(\boldsymbol{\beta}, \boldsymbol{\delta}_i) = {\bf Y}(\bbeta, \boldsymbol{\delta}_i)^T \mathcal{R}_{i}^*(\boldsymbol{\beta}, \boldsymbol{\delta}_i)^{-1}\bar{{\bf R}}_{i}(\boldsymbol{\beta}, \boldsymbol{\delta}_i)\mathcal{R}_{i}^*(\boldsymbol{\beta}, \boldsymbol{\delta}_i)^{-1}	 {\bf Y}(\bbeta, \boldsymbol{\delta}_i), \\	
{\bf M}_n^*(\boldsymbol{\beta}, \boldsymbol{\delta}) &:=& \sum_{i=1}^n {\bf Y}(\bbeta, \boldsymbol{\delta}_i)^T \bar{{\bf R}}_{i}(\boldsymbol{\beta}, \boldsymbol{\delta}_i)^{-1}	 {\bf Y}(\bbeta, \boldsymbol{\delta}_i), \\	
{\bf H}_n^*(\boldsymbol{\beta}, \boldsymbol{\delta}) &:=& \sum_{i=1}^n {\bf K}_i^*(\boldsymbol{\beta}, \boldsymbol{\delta}_i).\\
\end{eqnarray*}

Since ${\bf g}_n^s$ is optimal within $\mathcal{H}_n$, we show that
\begin{eqnarray}
&&\frac{\det {\bf H}_n^*(\boldsymbol{\beta})} {\det \bar{{\bf M}}_n(\boldsymbol{\beta})} \to 1, \quad \frac{\det {\bf H}_n^*(\boldsymbol{\beta}, \boldsymbol{\delta})} {\det \bar{{\bf M}}_n(\boldsymbol{\beta}, \boldsymbol{\delta})} \to 1 \quad \forall \boldsymbol{\beta} \in \mathcal{T}, \label{eq3} \\
&&\frac{\det {\bf M}_n^*(\boldsymbol{\beta})} {\det \bar{{\bf M}}_n(\boldsymbol{\beta})} \to 1, \quad \frac{\det {\bf M}_n^*(\boldsymbol{\beta}, \boldsymbol{\delta})} {\det \bar{{\bf M}}_n(\boldsymbol{\beta}, \boldsymbol{\delta})} \to 1 \quad \forall \boldsymbol{\beta} \in \mathcal{T}, \label{eq4}
\end{eqnarray}
which by Proposition 5.5 of \cite{heyde97} imply the conclusion.
We proceed with the proof of \eqref{eq3}; the proof of \eqref{eq4} is similar and is omitted. From $(O_4)$ we have,
\begin{equation}
\max \{\lambda_{\max}[\bar{{\bf L}}_i(\boldsymbol{\beta},\boldsymbol{\delta}_i)^{-1}], \lambda_{\max}[{\bf K}_i^*(\boldsymbol{\beta}, \boldsymbol{\delta}_i)^{-1}] \} < \infty \nonumber 
\end{equation} 
First, we show that 
\begin{equation}
\sup_{i \ge 1}\|\bar{{\bf L}}_i(\boldsymbol{\beta}, \boldsymbol{\delta}_i)^{-1}\| = \sup_{i \ge 1}\lambda_{\max}[\bar{{\bf L}}_i(\boldsymbol{\beta}, \boldsymbol{\delta}_i)^{-1}]< \infty, \label{norm.inv.L.bounded}
\end{equation}
or $\inf_{i \ge 1}\lambda_{\min}[\bar{{\bf L}}_i(\boldsymbol{\beta}, \boldsymbol{\delta}_i)] > 0.$ Indeed, let ${\bf x}$ be a vector of norm 1 such that ${\bf x}^T \bar{{\bf L}}_i(\boldsymbol{\beta}, \boldsymbol{\delta}_i){\bf x} =\lambda_{\min}[\bar{{\bf L}}_i(\boldsymbol{\beta}, \boldsymbol{\delta}_i)^{-1}].$ Since the absolute value of each entry of $\bar{{\bf R}}_i(\boldsymbol{\beta}, \boldsymbol{\delta}_i)$ is less than 1,
\begin{eqnarray*}
	\lambda_{\min}[\bar{{\bf L}}_i(\boldsymbol{\beta},\boldsymbol{\delta}_i)] &=& {\rm E}_{\boldsymbol{\beta}}[{\bf x}^T {\bf Y}_i(\boldsymbol{\beta}, \boldsymbol{\delta}_i)^T \bar{{\bf R}}_{i}(\boldsymbol{\beta}, \boldsymbol{\delta}_i)^{-1} {\bf Y}_i(\boldsymbol{\beta}, \boldsymbol{\delta}_i) {\bf x}] \\
	&\ge& {\rm E}_{\boldsymbol{\beta}}[\lambda_{\min}(\bar{{\bf R}}_{i}(\boldsymbol{\beta}, \boldsymbol{\delta}_i)^{-1}) {\bf x}^T {\bf Y}_i(\boldsymbol{\beta}, \boldsymbol{\delta}_i)^T {\bf Y}_i(\boldsymbol{\beta}, \boldsymbol{\delta}_i){\bf x}] \\
	&=& {\rm E}_{\boldsymbol{\beta}}\left[\frac{1}{\lambda_{\max}(\bar{{\bf R}}_{i}(\boldsymbol{\beta}, \boldsymbol{\delta}_i))} {\bf x}^T {\bf Y}_i(\boldsymbol{\beta}, \boldsymbol{\delta}_i)^T {\bf Y}_i(\boldsymbol{\beta}, \boldsymbol{\delta}_i){\bf x} \right] \\
	&\ge& \frac{1}{m} {\rm E}_{\boldsymbol{\beta}} \{ \lambda_{\min}[ {\bf Y}_i(\boldsymbol{\beta}, \boldsymbol{\delta}_i)^T {\bf Y}_i(\boldsymbol{\beta}, \boldsymbol{\delta}_i)]\}.
\end{eqnarray*}
Thus, $\inf_{i \ge 1} \lambda_{\min}[\bar{{\bf L}}_i(\boldsymbol{\beta},\boldsymbol{\delta}_i)]>0 $ by $(O_4).$

For any $p \times 1$ vector ${\bf x}$, we have the following inequalities
\begin{eqnarray}
&&\min_{n_0 \le i \le n} \lambda_{\min}[\bar{{\bf L}}_i(\boldsymbol{\beta},\boldsymbol{\delta}_i)^{-1}{\bf K}_i^*(\boldsymbol{\beta},\boldsymbol{\delta}_i)] {\bf x}^T\bar{{\bf M}}_{n_0,n}(\boldsymbol{\beta},\boldsymbol{\delta}){\bf x} \le {\bf x}^T{\bf H}_{n_0,n}^*(\boldsymbol{\beta},\boldsymbol{\delta}){\bf x} \nonumber\\
&&\le \max_{n_0 \le i \le n} \lambda_{\max}[\bar{{\bf L}}_i(\boldsymbol{\beta},\boldsymbol{\delta}_i)^{-1}{\bf K}_i^*(\boldsymbol{\beta},\boldsymbol{\delta}_i)] {\bf x}^T\bar{{\bf M}}_{n_0,n}(\boldsymbol{\beta},\boldsymbol{\delta}){\bf x}, \label{eval_M_H}
\end{eqnarray}
where $\bar{{\bf M}}_{n_0,n}(\boldsymbol{\beta},\boldsymbol{\delta}) := \sum_{i=n_0}^n \bar{{\bf L}}_i(\boldsymbol{\beta},\boldsymbol{\delta}_i)$ and ${\bf H}_{n_0,n}^*(\boldsymbol{\beta},\boldsymbol{\delta}) := \sum_{i=n_0}^n {\bf K}_i^*(\boldsymbol{\beta},\boldsymbol{\delta}_i).$
Next, we prove that as $ i \to \infty$,
\begin{equation}
\bar{{\bf L}}_i(\boldsymbol{\beta},\boldsymbol{\delta}_i)^{-1}{\bf K}_i^*(\boldsymbol{\beta},\boldsymbol{\delta}_i) \to {\bf I}. \label{convergenceLK}
\end{equation}
Combining \eqref{eval_M_H} and \eqref{convergenceLK}, we obtain that, for any $\varepsilon>0$, there exists $n_0= n_0(\varepsilon, \boldsymbol{\beta})$, such that for any $n_m \ge n_0$ and $n > n_m,$
\begin{equation}
(1-\varepsilon) {\bf x}^T\bar{{\bf M}}_{n_0,n}(\boldsymbol{\beta}, \boldsymbol{\delta}){\bf x} \le {\bf x}^T{\bf H}_{n_0,n}^*(\boldsymbol{\beta}, \boldsymbol{\delta}){\bf x}  \le (1+\varepsilon) {\bf x}^T\bar{{\bf M}}_{n_0,n}(\boldsymbol{\beta}, \boldsymbol{\delta}){\bf x}, \label{noua_eq}
\end{equation}
from which we derive
\begin{equation}
(1-\varepsilon)^p \le \frac{\det {\bf H}_{n_0,n}^*(\boldsymbol{\beta}, \boldsymbol{\delta})} {\det \bar{{\bf M}}_{n_0,n}(\boldsymbol{\beta}, \boldsymbol{\delta})} \le (1+\varepsilon)^p, \ \mbox{for any } n > n_m. \label{evaluation_partial_sum}
\end{equation}
We proceed with the proof of \eqref{convergenceLK}. We show that, as $i \to \infty$
\begin{equation}
\|\bar{{\bf L}}_i(\boldsymbol{\beta},\boldsymbol{\delta}_i)^{-1}{\bf K}_i^*(\boldsymbol{\beta},\boldsymbol{\delta}_i) - {\bf I}\| \le \|\bar{{\bf L}}_i(\boldsymbol{\beta},\boldsymbol{\delta}_i)^{-1}\| \cdot \|{\bf K}_i^*(\boldsymbol{\beta},\boldsymbol{\delta}_i)-\bar{{\bf L}}_i(\boldsymbol{\beta},\boldsymbol{\delta}_i)\| \to 0. \nonumber
\end{equation}
The first factor on the right hand side, is bounded, due to \eqref{norm.inv.L.bounded}. To complete the proof of \eqref{convergenceLK}, we now prove that $\|{\bf K}_i^*(\boldsymbol{\beta},\boldsymbol{\delta}_i)-\bar{{\bf L}}_i(\boldsymbol{\beta},\boldsymbol{\delta}_i)\| \to 0,$ $i \to \infty$.
We have
\begin{eqnarray*}
	&& \left\|{\rm E}_{\boldsymbol{\beta}}\left\{ {\bf Y}_i(\boldsymbol{\beta}, \boldsymbol{\delta}_i)^T \left[\mathcal{R}_{i}^*(\boldsymbol{\beta}, \boldsymbol{\delta}_i)^{-1}- \bar{{\bf R}}_{i}(\boldsymbol{\beta}, \boldsymbol{\delta}_i)^{-1}\right]{\bf Y}_i(\boldsymbol{\beta}, \boldsymbol{\delta}_i)\right\} \right\|  \\
	&& \le \left\|{\rm E}_{\boldsymbol{\beta}}\left\{ {\bf Y}_i(\boldsymbol{\beta}, \boldsymbol{\delta}_i)^T \left[\mathcal{R}_{i}^*(\boldsymbol{\beta}, \boldsymbol{\delta}_i)^{-1}- \mathcal{R}_{i}^*(\boldsymbol{\beta})^{-1}\right]{\bf Y}_i(\boldsymbol{\beta}, \boldsymbol{\delta}_i)\right\} \right\|  \\
	&& + \left\|{\rm E}_{\boldsymbol{\beta}}\left\{ {\bf Y}_i(\boldsymbol{\beta}, \boldsymbol{\delta}_i)^T \left[\mathcal{R}_{i}^*(\boldsymbol{\beta})^{-1}- \bar{{\bf R}}_{i}(\boldsymbol{\beta})^{-1} \right]{\bf Y}_i(\boldsymbol{\beta}, \boldsymbol{\delta}_i)\right\} \right\| \\
	&& + \left\|{\rm E}_{\boldsymbol{\beta}}\left\{ {\bf Y}_i(\boldsymbol{\beta}, \boldsymbol{\delta}_i)^T \left[\bar{R}_{i}^*(\boldsymbol{\beta})^{-1}- \bar{{\bf R}}_{i}(\boldsymbol{\beta}, \boldsymbol{\delta}_i)^{-1} \right]{\bf Y}_i(\boldsymbol{\beta}, \boldsymbol{\delta}_i)\right\} \right\|  \\
	&& := T_1(i) + T_2(i) + T_3(i). 
\end{eqnarray*}
Using $(O_2)$, we obtain the inequalities 
\begin{eqnarray*}
	&& T_1(i) \le 2{\rm E}_{\boldsymbol{\beta}} \left\{\left\| {\bf Y}_i(\boldsymbol{\beta}, \boldsymbol{\delta}_i)^T \left[\mathcal{R}_{i}^*(\boldsymbol{\beta}, \boldsymbol{\delta}_i)^{-1}- \mathcal{R}_{i}^*(\boldsymbol{\beta})^{-1}\right]{\bf Y}_i(\boldsymbol{\beta}, \boldsymbol{\delta}_i) \right\| \right\} \\
	&& \le 2{\rm E}_{\boldsymbol{\beta}} \left\{\| {\bf Y}_i(\boldsymbol{\beta}, \boldsymbol{\delta}_i) \|^2 \left\|\mathcal{R}_{i}^*(\boldsymbol{\beta}, \boldsymbol{\delta}_i)^{-1}- \mathcal{R}_{i}^*(\boldsymbol{\beta})^{-1}\right\| \right\} \le 2\frac{1}{2^i} \sup_{i \ge 1} {\rm E}_{\boldsymbol{\beta}} [\|{\bf Y}_i(\boldsymbol{\beta}, \boldsymbol{\delta}_i) \|^2].
\end{eqnarray*} 
By $(O_3),$ we obtain $T_1(i) \to 0,$ as $i \to \infty.$ On the other hand,  
\begin{eqnarray*}
	&& T_2(i) \le  2{\rm E}_{\boldsymbol{\beta}}\{\| {\bf Y}_i(\boldsymbol{\beta}, \boldsymbol{\delta}_i)^T [\mathcal{R}_{i}^*(\boldsymbol{\beta})^{-1}- \bar{{\bf R}}_{i}(\boldsymbol{\beta})^{-1} ]{\bf Y}_i(\boldsymbol{\beta}, \boldsymbol{\delta}_i) \|\} \\
	&& \le 2{\rm E}_{\boldsymbol{\beta}}\{\| {\bf Y}_i(\boldsymbol{\beta}, \boldsymbol{\delta}_i)\|^2 \|\mathcal{R}_{i}^*(\boldsymbol{\beta})^{-1}- \bar{{\bf R}}_{i}(\boldsymbol{\beta})^{-1} \|\} \\
	&& \le 2{\rm E}_{\boldsymbol{\beta}}\{\| {\bf Y}_i(\boldsymbol{\beta}, \boldsymbol{\delta}_i)\|^2 \|\mathcal{R}_{i}^*(\boldsymbol{\beta})^{-1}- \bar{{\bf R}}_{i}(\boldsymbol{\beta})^{-1} \| {\bf 1}_{\{\|\mathcal{R}_{i}^*(\boldsymbol{\beta})^{-1}- \bar{{\bf R}}_{i}(\boldsymbol{\beta})^{-1}\| < \varepsilon\} }\} \\
	&& +  2 {\rm E}_{\boldsymbol{\beta}}\{\| {\bf Y}_i(\boldsymbol{\beta}, \boldsymbol{\delta}_i)\|^2 \|\mathcal{R}_{i}^*(\boldsymbol{\beta})^{-1}- \bar{{\bf R}}_{i}(\boldsymbol{\beta})^{-1} \| {\bf 1}_{\{\|\mathcal{R}_{i}^*(\boldsymbol{\beta})^{-1}- \bar{{\bf R}}_{i}(\boldsymbol{\beta})^{-1}\| \ge \varepsilon\} }\}
\end{eqnarray*}
which gives
\begin{eqnarray*}
	T_2(i) \le 2\varepsilon \sup_{i \ge 1}{\rm E}_{\boldsymbol{\beta}}\{\| {\bf Y}_i(\boldsymbol{\beta}, \boldsymbol{\delta}_i)\|^2 \}+ 2M_1 \eta, \ i \ge i_0(\varepsilon).
\end{eqnarray*}
We obtained the bound for the second term on the right hand side by using conditions $(H'),$ $(R),$ $(O_1)$ and $(O_3).$ This proves that $T_2(i) \to 0$ as $i \to \infty.$ The proof of \eqref{convergenceLK} is concluded by remarking that 
\begin{eqnarray*}
	&& T_3(i) \le 2{\rm E}_{\boldsymbol{\beta}} \left\{\| {\bf Y}_i(\boldsymbol{\beta}, \boldsymbol{\delta}_i) \|^2 \left\|\bar{{\bf R}}_{i}(\boldsymbol{\beta})^{-1}- \bar{{\bf R}}_{i}(\boldsymbol{\beta}, \boldsymbol{\delta}_i)^{-1}\right\| \right\} \le 2 \varepsilon \sup_{i \ge 1} {\rm E}_{\boldsymbol{\beta}} [\|{\bf Y}_i(\boldsymbol{\beta}, \boldsymbol{\delta}_i) \|^2],
\end{eqnarray*} 
using the continuity of each entry of $\bar{{\bf R}}_{i}(\boldsymbol{\beta})$ (and hence, of its inverse) with respect to ${\bf X}_i$. 

Since all the eigenvalues of the matrix $\bar{{\bf L}}_i(\boldsymbol{\beta},\boldsymbol{\delta}_i)^{-1}{\bf K}_i^*(\boldsymbol{\beta},\boldsymbol{\delta}_i)$ converge to 1, for any $\varepsilon > 0,$ there exists $n_0$ such that 
\begin{equation}
1- \varepsilon \le \min_{n_0 \le i \le n} \lambda_{\min}[\bar{{\bf L}}_i(\boldsymbol{\beta},\boldsymbol{\delta}_i)^{-1}{\bf K}_i^*(\boldsymbol{\beta},\boldsymbol{\delta}_i)] \le \max_{n_0 \le i \le n} \lambda_{\max}[\bar{{\bf L}}_i(\boldsymbol{\beta},\boldsymbol{\delta}_i)^{-1}{\bf K}_i^*(\boldsymbol{\beta},\boldsymbol{\delta}_i)] \le 1+\varepsilon, \nonumber
\end{equation}
Combining these inequalities with \eqref{eval_M_H} we obtain \eqref{noua_eq}, and thus \eqref{evaluation_partial_sum}. Actually, our first goal is to obtain inequalities similar to those in \eqref{evaluation_partial_sum}, when all $\boldsymbol{\delta}_i=0,$ $i \ge 1.$ We proceed as follows.

For any $p \times 1$ vector ${\bf x}$, we have
\begin{eqnarray}
&&\min_{n_0 \le i \le n} \lambda_{\min}[{\bf K}_i^*(\boldsymbol{\beta}, \boldsymbol{\delta}_i)^{-1}{\bf K}_i^*(\boldsymbol{\beta})] {\bf x}^T{\bf H}_{n_0,n}^*(\boldsymbol{\beta}, \boldsymbol{\delta}){\bf x} \le {\bf x}^T{\bf H}_{n_0,n}^*(\boldsymbol{\beta}){\bf x} \nonumber\\
&&\le \max_{n_0 \le i \le n} \lambda_{\max}[{\bf K}_i^*(\boldsymbol{\beta}, \boldsymbol{\delta}_i)^{-1}{\bf K}_i^*(\boldsymbol{\beta})] {\bf x}^T{\bf H}_{n_0,n}^*(\boldsymbol{\beta}, \boldsymbol{\delta}){\bf x}. \label{eval_H_Hd}
\end{eqnarray}
and
\begin{eqnarray}
&&\min_{n_0 \le i \le n} \lambda_{\min}[\bar{{\bf L}}_i(\boldsymbol{\beta}, \boldsymbol{\delta}_i)^{-1}\bar{{\bf L}}_i(\boldsymbol{\beta})] {\bf x}^T\bar{{\bf M}}_{n_0,n}(\boldsymbol{\beta}, \boldsymbol{\delta}){\bf x} \le {\bf x}^T\bar{{\bf M}}_{n_0,n}(\boldsymbol{\beta}){\bf x} \nonumber\\
&&\le \max_{n_0 \le i \le n} \lambda_{\max}[\bar{{\bf L}}_i(\boldsymbol{\beta},\boldsymbol{\delta}_i)^{-1}\bar{{\bf L}}_i(\boldsymbol{\beta})] {\bf x}^T\bar{{\bf M}}_{n_0,n}(\boldsymbol{\beta}, \boldsymbol{\delta}){\bf x}. \label{eval_M_Md}
\end{eqnarray}
We will prove that, as $ i \to \infty$,
\begin{equation}
{\bf K}_i^*(\boldsymbol{\beta}, \boldsymbol{\delta}_i)^{-1}{\bf K}_i^*(\boldsymbol{\beta}) \longrightarrow {\bf I}, \label{convergenceKKd}
\end{equation}
and hence all the eigenvalues of the matrix ${\bf K}_i^*(\boldsymbol{\beta}, \boldsymbol{\delta}_i)^{-1}{\bf K}_i^*(\boldsymbol{\beta})$ converge to 1. If \eqref{convergenceKKd} holds then, for $\varepsilon>0$, there exists $n_1$ such that 
\begin{equation}
1- \varepsilon \le \min_{n_1 \le i \le n} \lambda_{\min}[{\bf K}_i^*(\boldsymbol{\beta}, \boldsymbol{\delta}_i)^{-1}{\bf K}_i^*(\boldsymbol{\beta})] \le \max_{n_1 \le i \le n} \lambda_{\max}[{\bf K}_i^*(\boldsymbol{\beta}, \boldsymbol{\delta}_i)^{-1}{\bf K}_i^*(\boldsymbol{\beta})] \le 1+\varepsilon. \nonumber
\end{equation}
These inequalities, combined with \eqref{eval_H_Hd} imply first that for all $n_m,$ $n,$ $n_1 \le n_m< n$
\begin{equation}
\label{eq222}
(1-\varepsilon) {\bf x}^T{\bf H}_{n_m,n}^*(\boldsymbol{\beta}, \boldsymbol{\delta}){\bf x} \le {\bf x}^T{\bf H}_{n_m,n}^*(\boldsymbol{\beta}){\bf x} \le (1+\varepsilon) {\bf x}^T{\bf H}_{n_m,n}^*(\boldsymbol{\beta}, \boldsymbol{\delta}){\bf x}
\end{equation}
and then 
\begin{equation}
\label{evaluation_partial_sum_H_Hd}
(1-\varepsilon)^p  \le \frac{\det {\bf H}_{n_m,n}^*(\boldsymbol{\beta})}{\det {\bf H}_{n_m,n}^*(\boldsymbol{\beta}, \boldsymbol{\delta})} \le (1+\varepsilon)^p, \ \mbox{for any } n > n_m. 
\end{equation}
Similarly, if 
\begin{equation}
\bar{{\bf L}}_i(\boldsymbol{\beta},\boldsymbol{\delta}_i)^{-1}\bar{{\bf L}}_i(\boldsymbol{\beta}) \longrightarrow {\bf I},  \label{new_label2}
\end{equation} 
then from \eqref{eval_M_Md} there exists $n_2$ such that, for any $n_m \ge n_2$ and $n \ge n_m,$ 
\begin{equation}
(1-\varepsilon){\bf x}^T\bar{{\bf M}}_{n_m,n}(\boldsymbol{\beta}, \boldsymbol{\delta}){\bf x} \le {\bf x}^T\bar{{\bf M}}_{n_m,n}(\boldsymbol{\beta}){\bf x} \le (1+\varepsilon){\bf x}^T\bar{{\bf M}}_{n_m,n}(\boldsymbol{\beta}, \boldsymbol{\delta}){\bf x}, \label{eq111}
\end{equation}
and 
\begin{equation}
\label{evaluation_partial_sum_M_Md}
(1-\varepsilon)^p \le \frac{ \det \bar{{\bf M}}_{n_m,n}(\boldsymbol{\beta})}{\det \bar{{\bf M}}_{n_m,n}(\boldsymbol{\beta}, \boldsymbol{\delta})} \le (1+\varepsilon)^p , \ \mbox{for any } n > n_m. 
\end{equation}
Let
\begin{equation}
\frac{\det {\bf H}_{n_m,n}^*(\boldsymbol{\beta})} {\det \bar{{\bf M}}_{n_m,n}(\boldsymbol{\beta})} = \frac{\det {\bf H}_{n_m,n}^*(\boldsymbol{\beta})} {\det {\bf H}_{n_m,n}^*(\boldsymbol{\beta}, \boldsymbol{\delta})} \frac{\det {\bf H}_{n_m,n}^*(\boldsymbol{\beta}, \boldsymbol{\delta})}{\det \bar{{\bf M}}_{n_m,n}(\boldsymbol{\beta}, \boldsymbol{\delta})} \frac{\det \bar{{\bf M}}_{n_m,n}(\boldsymbol{\beta}, \boldsymbol{\delta})}{\det \bar{{\bf M}}_{n_m,n}(\boldsymbol{\beta})}. \nonumber
\end{equation}
Combining \eqref{evaluation_partial_sum}, \eqref{evaluation_partial_sum_H_Hd} and \eqref{evaluation_partial_sum_M_Md} we obtain 
\begin{equation}
(1-\varepsilon)^p(1-\varepsilon)^{p}(1+\varepsilon)^{-p} \le \frac{\det {\bf H}_{n_m,n}^*(\boldsymbol{\beta})} {\det \bar{{\bf M}}_{n_m,n}(\boldsymbol{\beta})} \le (1+\varepsilon)^p(1+\varepsilon)^{p}(1-\varepsilon)^{-p}, \ \mbox{for any } n > n_m. \nonumber
\end{equation}
where $n_3 \le n_m \le n,$ where $n_3 = \max\{n_0, n_1, n_2\}.$
We now turn to the proof of \eqref{convergenceKKd} which closely follows the proof of \eqref{convergenceLK}. The proof of \eqref{new_label2} is similar. We show that, as $i \to \infty$
\begin{equation}
\|{\bf K}_i^*(\boldsymbol{\beta},\boldsymbol{\delta}_i)^{-1}{\bf K}_i^*(\boldsymbol{\beta}) - {\bf I}\| \le \|{\bf K}_i^*(\boldsymbol{\beta},\boldsymbol{\delta}_i)^{-1}\| \cdot \|{\bf K}_i^*(\boldsymbol{\beta})-{\bf K}_i^*(\boldsymbol{\beta},\boldsymbol{\delta}_i)\| \to 0. \nonumber
\end{equation}
By $(O_2)$ and $(O_3)$ we have $ \sup_{i \ge 1}{\rm E}_{\boldsymbol{\beta}}[ \|{\bf Y}_i(\boldsymbol{\beta})\|] < \infty.$\\

A further use of $(O_2)$ and $(O_3)$ gives
\begin{eqnarray*}
	&& \| {\bf K}_i^*(\boldsymbol{\beta})-{\bf K}_i^*(\boldsymbol{\beta},\boldsymbol{\delta}_i) \| \\
	&& = \| {\rm E}_{\boldsymbol{\beta}}[ {\bf Y}_i(\boldsymbol{\beta})^T \mathcal{R}_{i}^*(\boldsymbol{\beta})^{-1} {\bf Y}_i(\boldsymbol{\beta})] - {\rm E}_{\boldsymbol{\beta}}\left[ {\bf Y}_i(\boldsymbol{\beta}, \boldsymbol{\delta}_i)^T \mathcal{R}_{i}^*(\boldsymbol{\beta}, \boldsymbol{\delta}_i)^{-1} {\bf Y}_i(\boldsymbol{\beta}, \boldsymbol{\delta}_i)\right] \| \\
	&& \le 2{\rm E}_{\boldsymbol{\beta}} \|[{\bf Y}_i(\boldsymbol{\beta})^T - {\bf Y}_i(\boldsymbol{\beta}, \boldsymbol{\delta}_i)]^T\mathcal{R}_{i}^*(\boldsymbol{\beta})^{-1} {\bf Y}_i(\boldsymbol{\beta})\| \\
	&& +  2{\rm E}_{\boldsymbol{\beta}} \|{\bf Y}_i(\boldsymbol{\beta}, \boldsymbol{\delta}_i)^T [\mathcal{R}_{i}^*(\boldsymbol{\beta})^{-1}- \mathcal{R}_{i}^*(\boldsymbol{\beta}, \boldsymbol{\delta}_i)^{-1} ] {\bf Y}_i(\boldsymbol{\beta}, \boldsymbol{\delta}_i)\| \\
	&& + 2{\rm E}_{\boldsymbol{\beta}} \|{\bf Y}_i(\boldsymbol{\beta}, \boldsymbol{\delta}_i)^T \mathcal{R}_{i}^*(\boldsymbol{\beta})^{-1}[{{\bf A}_i}(\boldsymbol{\beta})^{1/2} {\bf X}_i -{\bf Y}_i(\boldsymbol{\beta}, \boldsymbol{\delta}_i) ] \| \\
	&\le& 2\frac{1}{2^i} K(\boldsymbol{\beta})^{-1} \sup_{i \ge 1}{\rm E}_{\boldsymbol{\beta}}[ \|{\bf Y}_i(\boldsymbol{\beta})\|] + 2\frac{1}{2^i} \sup_{i \ge 1} {\rm E}_{\boldsymbol{\beta}}[\|{\bf Y}_i(\boldsymbol{\beta}, \boldsymbol{\delta}_i)^2\|] + 4\frac{1}{2^i} K(\boldsymbol{\beta})^{-1} \sup_{i \ge 1}{\rm E}_{\boldsymbol{\beta}}[ \|{\bf Y}_i(\boldsymbol{\beta}, \boldsymbol{\delta}_i)\|]\\
	&& \to 0, \ \mbox{as }i \to \infty,
\end{eqnarray*}
by $(R)$ and because all expectations are equibounded. 
Thus \eqref{convergenceKKd} holds and the proof of \eqref{evaluation_partial_sum_H_Hd} is complete.\\ To complete the proof of the first part of \eqref{eq3}, we have to deal with the terms that are missing in \eqref{evaluation_partial_sum_H_Hd}.

Since $\bar{{\bf M}}_{n}(\boldsymbol{\beta}) \ge m {\bf H}_{n}^{\rm ind}(\boldsymbol{\beta}),$ $\lambda_{\min}[\bar{{\bf M}}_{n}(\boldsymbol{\beta})] \to \infty,$ as $n \to \infty$ by condition $(D^*),$ we obtain $\lambda_{\min}[\bar{{\bf M}}_{n_m, n}(\boldsymbol{\beta})] \to \infty,$ for any fixed $n_m$ as $n \to \infty.$ For a given $\varepsilon >0$ we can find $n_4 =n_4(\varepsilon,\boldsymbol{\beta}) > n_m,$ such that $\lambda_{\min}[\bar{{\bf M}}_{n_m, n}(\boldsymbol{\beta})] \ge \varepsilon^{-1} \lambda_{\max}[\bar{{\bf M}}_{n_m-1}(\boldsymbol{\beta})],$ for all $n \ge n_4.$ Thus, for all $n \ge n_4 > n_m.$

\begin{equation}
\nonumber
\bar{{\bf M}}_{n_m, n}(\boldsymbol{\beta}) \le \bar{{\bf M}}_{n}(\boldsymbol{\beta}) \le (1+\varepsilon)\bar{{\bf M}}_{n_m, n}(\boldsymbol{\beta}) 
\end{equation} 
and
\begin{equation}
\label{eq_new}
\det \bar{{\bf M}}_{n_m, n}(\boldsymbol{\beta})  \le  \det \bar{{\bf M}}_{n}(\boldsymbol{\beta})  \le  (1+\varepsilon)^p \det \bar{{\bf M}}_{n_m, n}(\boldsymbol{\beta}). 
\end{equation}
Combining \eqref{noua_eq}, \eqref{eq222} and \eqref{eq111} we obtain for all $n > n_m$
\begin{eqnarray*}
	(1+\varepsilon)^{-1} (1-\varepsilon)^2 {\bf x}^T\bar{{\bf M}}_{n_m, n}(\boldsymbol{\beta}){\bf x} \le {\bf x}^T{\bf H}_{n_m, n}^*(\boldsymbol{\beta}){\bf x} \le (1+\varepsilon)^2 (1-\varepsilon)^{-1} {\bf x}^T\bar{{\bf M}}_{n_m, n}(\boldsymbol{\beta}){\bf x}.
\end{eqnarray*}
These inequalities imply that $\lambda_{\min}[{\bf H}_{n_m,n}^*(\boldsymbol{\beta})] \to \infty.$ Reasoning as above, there exists an integer $n_5 \ge  n_4,$ such that, for all $n \ge n_5 > n_m,$
\begin{equation}
\det {\bf H}_{n_m, n}^*(\boldsymbol{\beta}) \le \det {\bf H}_{n}^*(\boldsymbol{\beta}) \le (1 + \varepsilon)^p \det {\bf H}_{n_m, n}^*(\boldsymbol{\beta}), \label{eq25}
\end{equation}
Combining \eqref{eq_new} and \eqref{eq25} gives, for $n \ge n_5 > n_m$
\begin{equation}
\frac{1} {(1+\varepsilon)^p} \frac{\det {\bf H}_{n_m, n}^*(\boldsymbol{\beta})}{\det \bar{{\bf M}}_{n_m,n}(\boldsymbol{\beta})}  \le \frac{\det {\bf H}_{n}^*(\boldsymbol{\beta})}{\det \bar{{\bf M}}_{n_m,n}(\boldsymbol{\beta})} \le (1+ \varepsilon)^p \frac{\det {\bf H}_{n_m, n}^*(\boldsymbol{\beta})}{\det \bar{{\bf M}}_{n_m,n}(\boldsymbol{\beta})}. \label{eq_eq}
\end{equation}
Finally, we obtain the first part of \eqref{eq3} from \eqref{evaluation_partial_sum_H_Hd} and \eqref{eq_eq}. To prove $\displaystyle{\frac{\det {\bf H}_{n}^*(\boldsymbol{\beta}, \boldsymbol{\delta})}{\det \bar{\bf M}_{n}(\boldsymbol{\beta}, \boldsymbol{\delta})}} \to 1,$ as $n \to \infty$ we proceed in a similar way. We note first that, by \eqref{eq222}, $\lambda_{\min}[{\bf H}_{m_n, n}^*(\boldsymbol{\beta})] \to \infty$ is equivalent to $\lambda_{\min}[{\bf H}_{m_n, n}^*(\boldsymbol{\beta}, \boldsymbol{\delta})] \to \infty,$ as $n \to \infty.$ Then $\bar{\bf M}_{m_n, n}(\boldsymbol{\beta})$ and ${\bf H}_{m_n, n}^*(\boldsymbol{\beta})$ can be replaced with $\bar{\bf M}_{m_n, n}(\boldsymbol{\beta}, \boldsymbol{\delta})$ and ${\bf H}_{m_n, n}^*(\boldsymbol{\beta}, \boldsymbol{\delta}),$ respectively in \eqref{eq_new} - \eqref{eq_eq}. Now, the second part of \eqref{eq3} follows from \eqref{evaluation_partial_sum}.  \hfill   $\Box$


\subsection{Verification of condition $(S)(i)$ of Theorem \ref{strong-consistency-th}}
\label{appen}
We define various moduli of continuity and then show that the asymptotic behaviour of combinations of these moduli are sufficient for $(S)(i)$ to hold.
For all $n \ge 1$, $r>0$, we define the following random variables
\begin{eqnarray*}
\gamma_n^{ '} = \max_{i \leq n, j \leq m_i} {\bf
x}_{ij}^T ({\bf H}_n')^{-1} {\bf x}_{ij}, \
a_n' = \lambda_{\rm max}({\bf H}_n') \gamma_n^{ '},
\end{eqnarray*}
\begin{eqnarray*}
\nu_n(r) &=& \sup_{\boldsymbol{\beta} \in B_r} \max_{i
	\leq n, j\le m_i} \left [\left|\frac{\mu^{''}({\bf x}_{ij}^T \boldsymbol{\beta}')}{\mu'({\bf x}_{ij}^T \boldsymbol{\beta})^{1/2}} - \frac{\mu^{''}({\bf x}_{ij}^T \boldsymbol{\beta}_0)}{\mu'({\bf x}_{ij}^T \boldsymbol{\beta}_0)^{1/2}}  \right| \cdot |\mu'({\bf x}_{ij}^T \boldsymbol{\beta}_0)|^{-1/2} \right], \\
\xi_n(r) &=&  \sup_{\boldsymbol{\beta} \in B_r} \max_{i
	\leq n, j\le m_i} \left [\left|\frac{\mu^{''}({\bf x}_{ij}^T \boldsymbol{\beta}')}{\mu'({\bf x}_{ij}^T \boldsymbol{\beta})^{3/2}} - \frac{\mu^{''}({\bf x}_{ij}^T \boldsymbol{\beta}_0)}{\mu'({\bf x}_{ij}^T \boldsymbol{\beta}_0)^{3/2}}  \right| \cdot |\mu'({\bf x}_{ij}^T \boldsymbol{\beta}_0)|^{1/2} \right], \\
 \pi_n(r) &=& \sup_{\boldsymbol{\beta} \in B_r} \max_{i
\leq n} \lambda_{\max}
[\cR_{i}^{1/2}\cR_{i}(\boldsymbol{\beta})^{-1}\cR_{i}^{1/2}],
 \\
 \rho_n(r) &=&  \sup_{\boldsymbol{\beta} \in B_r} \max_{i \leq n} \max_{j \le m_i}  \left| \lambda_j[\cR_{i}^{1/2}\cR_{i}(\boldsymbol{\beta})^{-1} \cR_{i}^{1/2} - \boldsymbol{I}] \right| \ \mbox{(note that } \rho_n(r) \mbox{ may converge to } 0,\\
 && \mbox{ as } r \to 0, n \to \infty),\\
d_n(r) &=&  \sup_{\boldsymbol{\beta} \in B_r}\max_{l \le p}\max_{i \le n} \max_{j \le m_i} \left|\lambda_{i}\left( \frac{\partial \cR_{i}(\boldsymbol{\beta})}{\partial \beta^l}\right) \right|, \\
\delta_n(r) &=&  \sup_{\boldsymbol{\beta} \in B_r}\max_{l \le p} \max_{i \le n} \max_{j \le m_i}\left|\lambda_{j} \left[ \left(\frac{\partial \cR_{i}(\boldsymbol{\beta}_0)}{\partial \beta^l}\right)^{-1/2} \left(\frac{\partial \cR_{i}(\boldsymbol{\beta})}{\partial \beta^l}\right) \left(\frac{\partial \cR_{i}(\boldsymbol{\beta}_0)}{\partial \beta^l}\right)^{-1/2} - {\bf I}\right] \right|.
\end{eqnarray*}
\noindent Note that $\delta_n(r) \mbox{ may converge to } 0, \mbox{ as }n \to \infty, \ r \to 0.$ \\
We remark that if 
\begin{equation}
(K) \hspace{7mm} \lim_{r \to 0} \limsup_{n \to \infty} r (a_n')^{1/2}=0, \nonumber
\end{equation}
then $\eta_n(r),$ $\nu_n(r)$ and $\xi_n(r)$ are bounded by $C r(a_n')^{1/2},$ for $r >r_0$. Furthermore, note also that $\pi_n(r) \ge 1,$ for all $n \ge 1$.\\

For $\displaystyle{0<\delta \le \frac{1}{2}}$, we introduce below a set of five conditions, which we label $(C)$.       
The strong conditions $(C_4)-(C_5)$ of Lemma 4.9 of \cite{balan-dumitrescu-schiopu10} are covered by the set of much weaker conditions stated below.
 \begin{eqnarray*}
(\gamma H') & & \limsup_{n \to \infty}(\gamma_n')^{1/2}[\lambda_{\max}({\bf H}_n')]^{1-\delta} < \infty  \ a.s., \\
(\pi) & & \lim_{r \to 0} \limsup_{n \to \infty} \pi_n(r) < \infty \  a.s., \\
(C_3') & & \lim_{r \to 0} \lim_{n \to \infty} r  d_n(r) [\lambda_{\rm max}({\bf H}_n')]^{1/2-\delta}=0 \ a.s., \\
(C_4) & & \lim_{r \to 0}\limsup_{n \to \infty} n E[r^2(a_n')^2  \pi_n(r)^2 \lambda_{\rm max}({\bf H}_n')]=0, \\
(C_5) & & \lim_{r \to 0}\limsup_{n \to \infty} n E[\widetilde{a}_n'   \rho_n(r)^2\lambda_{\rm max}({\bf H}_n')]=0, \ \widetilde{a}_n'=\max\{1, a_n'\},\\
(C_6) & & \lim_{r \to 0}\limsup_{n \to \infty} n E[r^2a_n' \pi_n(r)^4 d_n(r)^2 \lambda_{\rm max}({\bf H}_n')]=0, \\
(C_7) & & \lim_{r \to 0}\limsup_{n \to \infty} n E[\pi_n(r)^2 d_n(r)^2 \rho_n(r)^2\lambda_{\rm max}({\bf H}_n')]=0, \\
(C_8) & &\lim_{r \to 0} \limsup_{n \to \infty} n E[\pi_n(r)^2 \delta_n(r)^2\lambda_{\rm max}({\bf H}_n')]=0.
\end{eqnarray*} 

\begin{remark}
\label{NEWconditions}
{\rm Assume that $(\pi)$ and $(E')$ hold. Then $(C_1)$ and $(C_3)$ of \cite{balan-dumitrescu-schiopu10} simplify to $(\gamma H')$ and $(C_3')$, respectively. Furthermore, by Taylor's expansion, $(C_3')$ implies $(C_2)$ of \cite{balan-dumitrescu-schiopu10}, which can thus be replaced by the weaker condition $(\pi)$.}
\end{remark}

Next, we give conditions under which $(S)(i)$ of Theorem \ref{strong-consistency-th} holds.

\begin{prop}
If $(AH)$, $(E'),$ $(K)$ and $(C)$ hold, then $(S)(i)$ is satisfied, i.e.
\begin{equation}
\lim_{r \to 0} \limsup_{n \to \infty}
[\lambda_{\max}({\bf H}_n')]^{-1/2-\delta}\sup_{\boldsymbol{\beta} \in
B_r}\|| \cD_{n}(\boldsymbol{\beta})-\cD_{n}  \||=0, \ a.s. \nonumber
\end{equation}
\end{prop}
{\bf Proof.} The conclusion follows from the proofs of lemmas 4.6 - 4.9 in \cite{balan-dumitrescu-schiopu10}, under the new conditions, as shown in Remark \ref{NEWconditions}. Note that the normalizing factor is random, that condition $(K)$ in \cite{balan-dumitrescu-schiopu10} can be dropped and that their conditions $(E)$ and $(R')$ are satisfied here due to $(E')$. \hfill $\Box$

\vspace{1cm}

\begin{tabular}{l l}
\small{SCHOOL OF MATHEMATICS}                         & \small{DEPARTMENT OF MATHEMATICS  }\\
\ \small{AND STATISTICS}                               & \ \small{AND STATISTICS} \\
\small{VICTORIA UNIVERSITY OF WELLINGTON}              & \small{UNIVERSITY OF OTTAWA}\\
\small{WELLINGTON, 6140}                               & \small{OTTAWA, ONTARIO} \\
\small{NEW ZEALAND}           	                       & \small{CANADA K1N 6N5}\\
\end{tabular}

\end{document}